\documentclass[11pt]{article}
\usepackage{amsmath}
\usepackage{amsfonts}
\usepackage{color}
\usepackage{graphicx}
\usepackage{comment}
\usepackage{pgfplots}
\usepackage{caption}
\newcommand{\funding}[1]{Funding: #1}
\newenvironment{keywords}{\textbf{Key words:}}{}
\newenvironment{AMS}{\textbf{AMS:}}{}
\newenvironment{proof}{\textit{Proof.}}{}
\newtheorem{theo}{Theorem}[section]
\newtheorem{lem}[theo]{Lemma}
\newtheorem{prop}[theo]{Proposition}

\newtheorem{defi}{Definition}[section]

\newcommand{\mysection}[1]{\section{#1} \setcounter{equation}{0}}
\newcommand{\proofc}{{\sc Proof} \ }
\newcommand{\be}{\begin{equation} \label}
\newcommand{\ee}{\end{equation}}
\newcommand{\bea}{\begin{eqnarray}\label}
\newcommand{\eea}{\end{eqnarray}}
\newcommand{\bas}{\begin{eqnarray*}}
\newcommand{\eas}{\end{eqnarray*}}
\newcommand{\bit}{\begin{itemize}}
\newcommand{\eit}{\end{itemize}}
\newcommand{\qed}{\hfill$\Box$ \vskip.2cm}
\newcommand{\nn}{\nonumber}
\newcommand{\R}{\mathbb{R}}
\newcommand{\N}{\mathbb{N}}
\newcommand{\pO}{\partial\Omega}

\newcommand{\eps}{\varepsilon}

\newcommand{\wto}{\rightharpoonup}

\newcommand{\io}{\int_\Omega}
\newcommand{\bom}{\overline{\Omega}}
\newcommand{\abs}{\\[5pt]}

\newcommand{\ueps}{u_\eps}
\newcommand{\veps}{v_\eps}
\newcommand{\weps}{w_\eps}

\newcommand{\Feps}{F_\eps}

\newcommand{\F}{{\mathcal{F}}_\eps}
\newcommand{\D}{{\mathcal{D}}_\eps}
\newcommand{\hatu}{\widehat{U}}
\newcommand{\email}[1]{\texttt{#1}}

\begin{document}
\enlargethispage{10mm}
\title{Migration-driven benefit in a two-species nutrient taxis system}
\author{
	Piotr Krzy\.zanowski\thanks{Institute of Applied Mathematics and Mechanics, Warsaw University, ul.~Banacha 2, 02-097 Warszawa, Poland (\email{d.wrzosek@mimuw.edu.pl}, \email{p.krzyzanowski@mimuw.edu.pl})} %
	\and Michael Winkler\thanks{Institut f\"ur Mathematik, Universit\"at Paderborn, 33098 Paderborn, Germany (\email{michael.winkler@math.uni-paderborn.de}) \funding{The work of this author was supported by {\em Deutsche Forschungsgemeinschaft} in the context of the project {\em Analysis of chemotactic cross-diffusion in complex frameworks}}.} %
	\and Dariusz Wrzosek$^*$
	%\footnotemark[2] %
}
%\author{
%Michael Winkler\footnote{michael.winkler@math.uni-paderborn.de}\\
%{\small Institut f\"ur Mathematik, Universit\"at Paderborn,}\\
%{\small 33098 Paderborn, Germany} \\
%\medskip
%Piotr Krzy\.zanowski\footnote{P.Krzyzanowski@mimuw.edu.pl} and  Dariusz Wrzosek\footnote{ darekw@mimuw.edu.pl}\\
%{\small Institute of Applied Mathematics and Mechanics, Warsaw University,} \\
%{\small Banacha 2, 02-097 Warszawa, Poland} }
%%
%
\date{}
\maketitle
\begin{abstract}
A model describing the competition of two species for a common nutrient is studied. It is assumed that one of the
competitors is motionless while the other has the ability to move upwards gradients of the nutrient density. 
It is proved that under suitable assumptions on the initial data, in the long time perspective the ability to move 
turns out to be a crucial feature providing competitive advantage irrespectively of a possible difference between 
the species with regard to their rates of proliferation and nutrient uptake.
\end{abstract}

\mysection{Introduction}\label{intro}
We consider a basic  simplistic model describing the competition between two species feeding on a common single non-renewable resource -- a nutrient which is indispensable for reproduction albeit not necessary to survive. 
Our goal is to study possible benefits which may result from the ability to move in searching 
for food as a factor determining success in competition. 
Indeed we will see that in some cases this factor becomes crucial and the usual population parameters indicating 
competitive benefits become irrelevant. \abs
In particular, letting $w=w(x,t)$ denote the nutrient concentration and $u=u(x,t)$ and $v=v(x,t)$ represent the population 
densities of two competing species distributed  
in a bounded domain $\Omega\subset\R^n$, $n\ge 1\,,$ we shall consider the problem
\be{0A}
    	\left\{ \begin{array}{rcll}
	u_t &=& D_u \Delta u - \chi \nabla \cdot (u\nabla w) + \delta uw,
	\qquad & x\in\Omega, \ t>0, \\[1mm]
	v_t &=& \alpha vw,
	\qquad & x\in\Omega, \ t>0, \\[1mm]
	w_t &=& D_w\Delta w - \beta uw - \gamma vw,
	\qquad & x\in\Omega, \ t>0, \\[1mm]
	& & \hspace*{-14mm}
	\frac{\partial u}{\partial\nu}=\frac{\partial w}{\partial\nu}=0,
	\qquad & x\in\pO, \ t>0, \\[1mm]
	& & \hspace*{-14mm}
	u(x,0)=u_0(x), \quad
	v(x,0)=v_0(x), \quad
	w(x,0)=w_0(x),
	\qquad & x\in\Omega,
	\end{array} \right.
\ee
where $D_u>0$ and $D_w>0$ are diffusion constants and $\chi>0$ is the chemotactic sensitivity coefficient.  
The coefficient $\beta$  describes the consumption rate of the first species  while 
$\delta$  determines the species specific rate of proliferation and similar interpretation applies  to  the second species. \abs
Focusing on the interplay between population growth on one hand, and mobility enhancing effectiveness of foraging on the
other hand we take into account neither mortality nor nutrient renewal which are usually considered in the class of models
describing the competition for common resources. We underline that neglecting mortality terms in the model can be 
justified through a restriction to not too long time scales in the case of bacteria which after all proliferate by cell 
division. 
It should moreover be noted that in some related systems studied in the literature the rate of consumption is modelled by
Monod-type functions accounting for effects of saturation which inevitably occur at higher food densities; 
as in our case the nutrient is being continuously depleted it will never reach higher levels of density at which any such
saturation effect could become significant. For simplicity in presentation, we may therefore we restrict to 
bilinear Lotka-Volterra interaction noticing that our analysis can be extended to cases of widely arbitrary, and in particular
to bounded, functional consumption rates.\abs
The ODE part of the above system refers to model B studied in \cite{Hsu1} in which logistic growth of the nutrient and 
mortality of competing species are assumed to be  negligible  in the time scale of our interest where the impact  of species motility on feeding effectiveness is taken into account. The case of renewable resources and mortality in competing species demands a  separate study. 
We refer a reader to Tilman's monograph \cite{Tilman} which contain basic ODE  models of species competition for common 
resources particularly being  related to  experiments in a flow reactor named chemostat, also known as bio-reactor.\abs
The classical ODE model for chemostat dynamics with constant supply of nutrients and stirring was analyzed in \cite{Hsu1} 
where the long time behavior was described and especially competitive exclusion shown. 
Next, the corresponding reaction-diffusion system additionally accounting for random diffusion in an unstirred chemostat 
was studied e.g.~in \cite{Hsu2} and \cite{DL0}. The role of chemotaxis in the process of nutrient competition was pointed out
in the pioneering work \cite{DKL}, and then the resulting competition system with nutrient taxis in the 
spatially one-dimensional case in \cite{XWang} and in \cite{DL2}, where inter alia a global attractor 
has been shown to exist.\abs
It can be easily proved (see Appendix A) that in the case of the migration-free ODE system associated with (\ref{0A}), 
and hence also for the PDE system (\ref{0A}) with
constant initial data, each individual trajectory $(u,v,w)$ approaches an equilibrium of the form
$(u_\infty,v_\infty,0)$ with some nonnegative constants $u_\infty$ and $v_\infty$;
moreover, if the initial densities of both competitors are equal in that $u_0=v_0$,
then the population with the higher proliferation coefficient asymptotically dominates in the sense that
\be{0win}
	{\rm{sgn}}(u_\infty-v_\infty) ={\rm{sgn}}(\delta -\alpha)\,. 
\ee
The main intention of this work is to provide some analytical  evidence indicating a substantial change of this
picture when one of the competitors has the ability to move randomly and chemotactically 
towards increasing nutrient concentrations, whereas the other one remains sessile. 
This phenomenon 
is illustrated by numerical simulations which  show the  complexity of possible solutions behaviour depending on 
initial data and model parameters. 
Motility-driven beneficial effects in competitive biological systems
were detected experimentally e.g.~in the case of bacterial soilborne 
plant pathogen (\cite{yao}, \cite{Julie}); 
however, we are not aware of any analytical study rigorously confirming the occurrence of phenomena of this form
in (\ref{0A}) or any related model.\abs
{\bf Main results.} \quad
In order to formulate our main results in this direction in a conveniently simple framework, 
let us firstly rescale (\ref{0A}) by noting that letting
$\tilde{x}=\frac{x}{L}, \tilde{t}=\frac{t}{\tau}, u=U\tilde{u}, v=V\tilde{v}$ and $w=W\tilde{w}$
%\[
%	\tilde{x}=\frac{x}{L}\,, \; \tilde{t}=\frac{t}{\tau}\,,\;u=U\tilde{u}\,,\,,\;v=V\tilde{v}\,,\;w=W\tilde{w}
%\]
as well as
$D=\frac{D_u \tau}{L^2}, 1=D^\prime_w=\frac{D_w W\tau}{L^2}=\delta^\prime=\delta\tau W,
\chi^\prime=\frac{\chi U W\tau}{L^2}, \alpha^\prime =\alpha^\prime \tau  W,
\beta^\prime=\beta \tau U$ and $\gamma^\prime=\gamma \tau U$
%\[
%	D=\frac{D_u \tau}{L^2}\,,\;1=D^\prime_w=\frac{D_w W\tau}{L^2}=\delta^\prime=\delta\tau W \,, 
%	\chi^\prime=\frac{\chi U W\tau}{L^2}\,, \alpha^\prime =\alpha^\prime \tau  W\,,\; 
%	\beta^\prime=\beta \tau U \,,\;\gamma^\prime=\gamma \tau U 
%\]
for suitably chosen positive numbers $L$, $\tau$, $U,V$ and $W$,
on performing standard computations and dropping tildes and primes we arrive at the normalized version of (\ref{0}) given
by
\be{0}
    	\left\{ \begin{array}{rcll}
	u_t &=& D \Delta u - \chi \nabla \cdot (u\nabla w) + uw,
	\qquad & x\in\Omega, \ t>0, \\[1mm]
	v_t &=& \alpha vw,
	\qquad & x\in\Omega, \ t>0, \\[1mm]
	w_t &=& \Delta w - \beta uw - \gamma vw,
	\qquad & x\in\Omega, \ t>0, \\[1mm]
	& & \hspace*{-14mm}
	\frac{\partial u}{\partial\nu}=\frac{\partial w}{\partial\nu}=0,
	\qquad & x\in\pO, \ t>0, \\[1mm]
	& & \hspace*{-14mm}
	u(x,0)=u_0(x), \quad
	v(x,0)=v_0(x), \quad
	w(x,0)=w_0(x),
	\qquad & x\in\Omega,
	\end{array} \right.
\ee
which will be the particular objective of our subsequent considerations.
Here for the initial data, throughout the sequel we will suppose that their first two components are such that
\be{iuv}
    	\left\{	\begin{array}{l}
	u_0\in C^0 \mbox{ with $u_0>0$ in $\bom$}
	\quad \mbox{and} \\[1mm]
	v_0\in W^{2,\infty}(\Omega) \mbox{ satisfies $v_0>0$ in } \bom,
    	\end{array} \right.
\ee
and that for the signal density we have
\be{iw}
	w_0 \in W^{1,\infty}(\Omega) \mbox{ with $w_0>0$ in } \bom.
\ee
Constituting an apparently necessary prerequisite for any qualitative analysis,
the first of our results then asserts global solvability within a natural weak solution concept.
\begin{theo}\label{theo222}
  Let $n\le 5$ and $\Omega\subset\R^n$ be a bounded convex domain with smooth boundary, and suppose that 
  $D, \chi, \alpha, \beta$ and $\gamma$ are positive.
  Then for any choice of $u_0, v_0$ and $w_0$ fulfilling (\ref{iuv}) and (\ref{iw}),
  there exist nonnegative functions
  \be{222.1}
	\left\{ \begin{array}{l}
	u\in L^\frac{n+2}{n}_{loc}(\bom\times [0,\infty)) 
	\cap L^\frac{n+2}{n+1}_{loc}([0,\infty);W^{1,\frac{n+2}{n+1}}(\Omega)), \\[1mm]
	v \in L^\infty(\Omega\times (0,\infty)) \qquad \mbox{and} \\[1mm]
	w\in L^4_{loc}([0,\infty); W^{1,4}(\Omega))
	\end{array} \right.
  \ee
  such that $(u,v,w)$ is a global weak solution of (\ref{0}) in the sense of Definition \ref{defi_weak}.
\end{theo}
We remark that in the case $n\le 2$,
by means of an adequate extension of the argument given in \cite{taowin_consumption} (cf.~also
the outline in \cite[Section 7]{win_ct_abs}) it is possible to show that the {\em a priori} estimates derived in
Section \ref{sect2} below are actually sufficient to ensure global solvability even in the context of classical 
solutions. Since our focus will rather be on a qualitative analysis of (\ref{0}), however, we refrain from
pursuing this any further here.\abs
%
%
%
%
%Now approaching the main outcome of this study,
%for the following result on qualitative behavior of the above solutions, in addition to (\ref{iw})
%we introduce the further requirements on $w_0$ given by
%\be{M}
%	\io \frac{|\nabla w_0|^2}{w_0} \le M
%\ee
%and
%\be{d}
%	\|w_0\|_{L^\infty(\Omega)} \le \delta
%\ee
%for $M>0$ and $\delta>0$.
%
%
Now the core of this study, to be addressed in Section \ref{sect3},
reveals that relative to the second subpopulation, the first indeed 
may eventually take advantage of its ability to migrate, as becoming manifest in our main result on qualitative behavior
in (\ref{0}):
\begin{theo}\label{theo18}
  Let $n\le 5$, $D>0$, $\chi>0$,
  $\alpha>0, \beta>0$ and $\gamma>0$, and suppose that $u_0$ and $v_0$ are such that (\ref{iuv}) holds
  with $u_0\not\equiv const.$. 
  Then there exist $C>0$ and $T>0$ such that for all $M>0$ one can find $\sigma>0$ with the property that whenever
  $w_0$ complies with (\ref{iw}) and moreover satisfies
  \be{M}
	\io \frac{|\nabla w_0|^2}{w_0} \le M
  \ee
  and
  \be{d}
	\|w_0\|_{L^\infty(\Omega)} \le \sigma,
  \ee
% (\ref{M}) and (\ref{d}), 
  for the global weak solution $(u,v,w)$ of (\ref{0}) from
  Theorem \ref{theo222} we have $\ln \frac{v(\cdot,t)}{u(\cdot,t)} \in L^1(\Omega)$ for a.e.~$t>T$ and
  \be{18.1}
	I(t):=\io \ln \frac{v(x,t)}{u(x,t)}\, dx \le \io \ln \frac{v_0}{u_0} - C
	\qquad \mbox{for a.e. } t>T.
  \ee
\end{theo}
Particularly, in the special case of precisely identical initial densities $u_0$ and $v_0$, not necessarily spatially constant,
(\ref{18.1}) states that for suitably small initial nutrient distributions, the first population 
ultimately prevails in the sense that then the corresponding logarithmic averages satisfy
$\frac{1}{|\Omega|} \io \ln u \ge \frac{1}{|\Omega|} \io \ln v + C$ for all sufficiently large times, with a fixed
positive constant $C$.
We emphasize that in stark contrast to the migration-free asymptotics characterized by (\ref{0win}), through
its qualitative independence of the sizes of $\alpha$,  $\beta$ and $\gamma$
this result especially covers arbitrarily large values of $\alpha$ and hence applies to situations in which
the static species, thus equipped with a substantial advantage in the efficiency of proliferation,
might be expected to retain certain benefits thereof.\\
Section 4 finally contains results of numerical simulations performed in the case $n=1$, and for some radially symmetric 
solutions in three-dimensional balls.
Going partially beyond situations addressed in Theorem \ref{theo18}, these findings will indicate that in general 
situations not necessarily complying with the hypotheses from Theorem \ref{theo18}, 
the final result of competition may heavily depend on a subtle interplay of many factors, including the initial distribution of species and nutrient, while the increase of space dimension seems to affect formation of more spiky 
space-time patterns.		% indicated by the increase of the solutions gradient. 
\mysection{Regularization, a quasi-energy property and global solvability}\label{sect2}
\subsection{A weak solution concept and a family of approximate problems}
In our analysis we shall pursue the following rather natural concept generalizing the notion of solution to (\ref{0}).
\begin{defi}\label{defi_weak}
  Suppose that
  $u\in L^1_{loc}([0,\infty);W^{1,1}(\Omega)),
  v\in L^1_{loc}(\bom\times [0,\infty))$ and 
  $w\in L^\infty_{loc}(\bom\times [0,\infty)) \cap L^1_{loc}([0,\infty);W^{1,1}(\Omega))$
%  \be{w1}
%	\left\{ \begin{array}{l}
%	u\in L^1_{loc}([0,\infty);W^{1,1}(\Omega)), \\[1mm]
%	v\in L^1_{loc}(\bom\times [0,\infty))
%	\qquad \mbox{and} \\[1mm]
%	w\in L^\infty_{loc}(\bom\times [0,\infty)) \cap L^1_{loc}([0,\infty);W^{1,1}(\Omega))
%	\end{array} \right.
%  \ee
  are all nonnegative and such that
  $u\nabla w\in L^1_{loc}(\bom\times [0,\infty);\R^n)$.
%  \be{w2}
%	u\nabla w\in L^1_{loc}(\bom\times [0,\infty);\R^n).
%  \ee
  Then $(u,v,w)$ will be called a {\em global weak solution} of (\ref{0}) if for all
  $\varphi\in C_0^\infty(\bom\times [0,\infty))$, the identities
  \be{w3}
	- \int_0^\infty \io u\varphi_t
	- \io u_0\varphi(\cdot,0)
	= - D \int_0^\infty \io \nabla u \cdot \nabla \varphi
	+ \chi \int_0^\infty \io u\nabla w\cdot\nabla \varphi
	+ \int_0^\infty \io uw\varphi
  \ee
  and
  \be{w4}
	- \int_0^\infty \io v\varphi_t - \io v_0\varphi(\cdot,0)
	=\alpha \int_0^\infty \io vw\varphi
  \ee
  as well as
  \be{w5}
	-\int_0^\infty \io w\varphi_t - \io w_0\varphi(\cdot,0)
	= - \int_0^\infty \io \nabla w\cdot\nabla\varphi
	- \beta \int_0^\infty \io uw\varphi
	- \gamma \int_0^\infty \io vw\varphi
  \ee
  are satisfied.
\end{defi}
In order to construct such a global weak solutions by means of a suitable approximation procedure,
following \cite{win_ct_abs} we fix a family $(\Feps)_{\eps\in (0,1)} \subset C^\infty([0,\infty))$ 
of functions such that whenever $\eps\in (0,1)$, we have
\be{F1}
	\Feps(0)=0
	\qquad \mbox{and} \qquad
	0 \le \Feps'(s) \le 1
	\quad \mbox{for all } s\ge 0,
\ee
that
$\Feps \mbox{ and } 0\le s\mapsto s\Feps'(s)
\mbox{ are bounded on $[0,\infty)$ for each $\eps\in (0,1)$, with }
\Feps'(s) \nearrow 1
\mbox{ as } \eps\searrow 0$
for all $s\ge 0$,
%\be{F2}
%	\left\{
%	\begin{array}{l}
%	\Feps \mbox{ and } 0\le s\mapsto s\Feps'(s)
%	\quad \mbox{ are bounded on $[0,\infty)$ for each $\eps\in (0,1)$, with} \nn\\[2mm]
%	\Feps'(s) \nearrow 1
%	\quad \mbox{as } \eps\searrow 0
%	\quad \mbox{for all } s\ge 0,
%	\end{array}  \right.
%\ee
noting that these requirements are met if e.g.~we define
\bas
	\Feps(s):=\frac{s}{1+\eps s}
	\qquad \mbox{for $\eps\in (0,1)$ and } s\ge 0.
\eas
For $\eps\in (0,1)$, we then consider 
\be{0eps}
    	\left\{ \begin{array}{rcll}
	u_{\eps t} &=& D \Delta \ueps - \chi \nabla \cdot \Big( \ueps \Feps'(\ueps)\nabla \weps \Big) 
		+ \Feps(\ueps)\weps,
	\qquad & x\in\Omega, \ t>0, \\[1mm]
	v_{\eps t} &=& \alpha \veps \weps, 
	\qquad & x\in\Omega, \ t>0, \\[1mm]
	w_{\eps t} &=& \Delta \weps - \beta \Feps(\ueps) \weps -\gamma \veps\weps,
	\qquad & x\in\Omega, \ t>0, \\[1mm]
	& & \hspace*{-15mm} 
	\frac{\partial\ueps}{\partial\nu}=\frac{\partial\weps}{\partial\nu}=0,
	\qquad & x\in\pO, \ t>0, \\[1mm]
	& & \hspace*{-15mm} 
	\ueps(x,0)=u_0(x), \quad
	\veps(x,0)=v_0(x), \quad
	\weps(x,0)=w_0(x),
	\qquad & x\in\Omega,
	\end{array} \right.
\ee
By adapting essentially well-known arguments, as detailed e.g.~in \cite{taowin_PROCA} for a related problem and thus
omitted here,
one can readily verify that all these problems admit globally defined classical solutions:
\begin{lem}\label{lem_loc}
  Assume (\ref{iuv}) and (\ref{iw}).
  Then for each $\eps\in (0,1)$, there exist functions
  \be{l1}
	\left\{ \begin{array}{l}
	\ueps \in C^0(\bom\times [0,\infty)) \cap C^{2,1}(\bom\times (0,\infty)), \\[1mm]
	\veps \in C^1(\bom\times [0,\infty)), \\[1mm]
	\weps \in \bigcap_{p>n} C^0([0,\infty);W^{1,p}(\Omega)) \cap C^{2,1}(\bom\times (0,\infty)),
	\end{array} \right.
  \ee
  all positive in $\bom\times [0,\infty)$,
  such that $(\ueps,\veps,\weps)$ solves (\ref{0eps}) in the classical sense
  in $\Omega\times (0,\infty)$.
\end{lem}
The following bound for the total mass in the first solution component is elementary but important.
\begin{lem}\label{lem02}
  Assume (\ref{iuv}) and (\ref{iw}). Then for all $\eps\in (0,1)$,
  \be{mass}
	\io \ueps(\cdot,t) \le \io u_0 + \frac{1}{\beta} \io w_0
	\qquad \mbox{for all } t>0.
  \ee
\end{lem}
\proof
  We integrate the first and the third equation from (\ref{0eps}) over $x\in\Omega$ to see that
  \bas
	\frac{d}{dt} \io \ueps = \io \Feps(\ueps) \weps
	\qquad \mbox{for all } t>0
  \eas
  and
  \bas
	\frac{d}{dt} \io \weps 
	= - \beta \io \Feps(\ueps)\weps - \gamma \io \veps\weps
	\le -\beta\io \Feps(\ueps)\weps
	\qquad \mbox{for all } t>0.
  \eas
  Taking an appropriate linear combination and integrating in time shows that
  \bas
	\io \ueps(\cdot,t) + \frac{1}{\beta} \io \weps(\cdot,t)
	\le \io u_0 + \frac{1}{\beta} \io w_0
	\qquad \mbox{for all } t>0
  \eas
  and thus implies (\ref{mass}).
\qed
In view of the comparatively simple structures of the second and third equations in (\ref{0eps}), for the respective
solution components some basic information can even be obtained in a pointwise sense:
\begin{lem}\label{lem2}
  Suppose that (\ref{iuv}) and (\ref{iw}) hold. Then for all $\eps\in (0,1)$,
  \be{vinf}
	\inf_{y\in\Omega} v_0(y) \le \veps(x,t) 
	\le \|v_0\|_{L^\infty(\Omega)} \cdot e^{\frac{\alpha}{\kappa} \|w_0\|_{L^\infty(\Omega)}}
	\qquad \mbox{for all $x\in\Omega$ and } t>0
  \ee
  and
  \be{winfty}
	\|\weps(\cdot,t)\|_{L^\infty(\Omega)}
	\le \|w_0\|_{L^\infty(\Omega)} e^{-\kappa t}
	\qquad \mbox{for all } t>0,
  \ee
  where
  \be{kappa}
	\kappa:=\gamma \cdot \inf_{y\in\Omega} v_0(y) >0.
  \ee
\end{lem}
\proof
  Since $v_{\eps t}\ge 0$ by (\ref{0eps}), the left inequality in (\ref{vinf}) is obvious.
  As therefore
  \bas
	w_{\eps t} \le \Delta \weps - \gamma \veps\weps
	\le \Delta\weps - \kappa\weps
	\qquad \mbox{in } \Omega\times (0,\infty)
  \eas
  due to (\ref{F1}), we next obtain (\ref{winfty}) from the maximum principle.
  Since $\kappa$ is positive, on integration of the second equation in (\ref{0eps}) this in turn implies that
  \bas
	\veps(x,t)
	&=& v_0(x) \cdot e^{\alpha \int_0^t \weps(x,s) ds} \\
	&\le& \|v_0\|_{L^\infty(\Omega)} \cdot \exp \bigg\{ \alpha \int_0^t \|w_0\|_{L^\infty(\Omega)} e^{-\kappa s} ds
		\bigg\} \\
	&\le& \|v_0\|_{L^\infty(\Omega)} \cdot \exp \Big\{ \alpha \|w_0\|_{L^\infty(\Omega)} \cdot
		\frac{1-e^{-\kappa t}}{\kappa} \Big\}
	\qquad \mbox{for all $x\in\Omega$ and } t>0
  \eas
  and thus establishes also the right inequality in (\ref{vinf}).
\qed
\subsection{Constructing a quasi-energy functional}
%\subsection{Global existence in (\ref{0}). Proof of Theorem \ref{theo222}}
%
%
%
%
%
%
%
%
Now the core of both our existence theory as well as our subsequent qualitative analysis will be formed by the following
observation on presence of a favorable global quasi-dissipative structure in (\ref{0eps}). 
\begin{lem}\label{lem4}
  Let $K>0$. Then there exists $C(K)>0$ such that if (\ref{iuv}) and (\ref{iw}) hold as well as
  \be{4.1}
	\|w_0\|_{L^\infty(\Omega)} \le K,
  \ee
  then for each $\eps\in (0,1)$,
  \be{4.2}
	\F(t) := \beta \io \ueps\ln \ueps
	+ \frac{\gamma\chi}{2\alpha} \io \frac{|\nabla \veps|^2}{\veps}
	+ \frac{\chi}{2} \io \frac{|\nabla \weps|^2}{\weps}, 
	\qquad t\ge 0,
  \ee
  and
  \be{4.3}
	\D(t) := \io \frac{|\nabla \ueps|^2}{\ueps}
	+ \io |\Delta \weps|^2 
 	+ \io |\nabla \weps|^4,
	\qquad t>0,
  \ee
  have the property that for any $t>0$ we have
  \bea{4.4}
	\frac{d}{dt} \F	%(t)\Big[\ueps(\cdot,t),\veps(\cdot,t),\weps(\cdot,t)\Big]
	+ \frac{1}{C(K)} \cdot \D\Big[\ueps(\cdot,t),\weps(\cdot,t)\Big] 
	\le C(K) e^{-\kappa t} \cdot \Big\{ 1+\F\Big[\ueps(\cdot,t),\veps(\cdot,t),\weps(\cdot,t)\Big]\Big\},
  \eea
  where $\kappa>0$ is taken from (\ref{kappa}).
\end{lem}
\proof
  By using (\ref{0eps}), we compute
  \bea{4.5}
	\frac{d}{dt} \io \ueps\ln\ueps
%	= \io u_{\eps t} \ln \ueps + \io u_{\eps t} \nn\\
	= - D \io \frac{|\nabla\ueps|^2}{\ueps} 
	+ \chi \io \Feps'(\ueps) \nabla\ueps\cdot\nabla\weps
	+ \io \Feps(\ueps) \ln\ueps \cdot \weps  
	+ \io \Feps(\ueps)\weps
  \eea
  and
  \bea{4.6}
	\frac{d}{dt} \io \frac{|\nabla\veps|^2}{\veps}
%	&=& 2\io \frac{1}{\veps} \nabla\veps \cdot \nabla (\alpha\veps\weps)
%	- \io \frac{|\nabla\veps|^2}{\veps^2} \cdot\alpha\veps\weps \nn\\
	= \alpha\io \frac{\weps}{\veps} |\nabla\veps|^2
	+ 2\alpha \io \nabla\veps\cdot\nabla\weps
%	\qquad \mbox{for all } t>0
  \eea
  as well as
  \bea{4.7}
	\frac{d}{dt} \io \frac{|\nabla\weps|^2}{\weps}
%	&=& 2\io \frac{1}{\weps} \nabla\weps \cdot 
%	\nabla \Big(\Delta\weps - \beta \Feps(\ueps) \weps - \gamma\veps\weps\Big) \nn\\
%	& & 
%	- \io \frac{|\nabla\weps|^2}{\weps^2} \cdot \Big( \Delta\weps - \beta\Feps(\ueps) \weps - \gamma\veps\weps\Big)
%		\nn\\
	&=& - \io \frac{1}{\weps} \Delta |\nabla\weps|^2
	- 2\io \frac{1}{\weps} |D^2\weps|^2
	- \io \frac{|\nabla\weps|^2}{\weps^2} \Delta\weps \nn\\
	& & - \beta \io \frac{\Feps(\ueps)}{\weps} |\nabla\weps|^2
	- 2\beta \io \Feps'(\ueps) \nabla\ueps\cdot\nabla\weps \nn\\
	& & - \gamma \io \frac{\veps}{\weps} |\nabla\weps|^2
	-2\gamma \io \nabla\veps\cdot\nabla\weps
%	\qquad \mbox{for all } t>0.
  \eea
  for all $t>0$.
  Here we recall that for all positive $\varphi\in C^3(\bom)$ fulfilling $\frac{\partial\varphi}{\partial\nu}=0$
  on $\pO$, by straightforward computation relying on the fact that then
  $\frac{\partial |\nabla\varphi|^2}{\partial\nu}\le 0$ on $\pO$ by convexity of
  $\Omega$ (\cite{lions_ARMA}),
  \bas
%	& & \hspace*{-20mm}
	- \io \frac{1}{\varphi} \Delta |\nabla\varphi|^2
	- 2\io \frac{1}{\varphi} |D^2\varphi|^2
	- \io \frac{|\nabla\varphi|^2}{\varphi^2} \Delta\varphi
%	\le - \io \frac{1}{\varphi^2} \nabla\varphi \cdot \nabla |\nabla\varphi|^2
%	- 2\io \frac{1}{\varphi} |D^2\varphi|^2
%	- \io \frac{|\nabla\varphi|^2}{\varphi^2} \Delta\varphi \\
	\le - 2 \io \weps |D^2\ln \varphi|^2,
%	\qquad \mbox{for all } t>0,
  \eas
  and that there exist $c_1>0$ and $c_2>0$ with the property that for any such $\varphi$ we have
  \bas
	\io \frac{|\Delta\varphi|^2}{\varphi} \le c_1 \io \varphi |D^2\ln \varphi|^2
  \eas
  as well as
  \bas
	\io \frac{|\nabla\varphi|^4}{\varphi^3}
	\le c_2 \io \varphi |D^2\ln\varphi|^2
  \eas
  (cf.~\cite[Section 3]{win_CPDE2012}).
  Therefore, combining (\ref{4.5})-(\ref{4.7}) shows that 	%with $\F$ as in (\ref{4.2}) we have
  \bea{4.8}
	& & \hspace*{-20mm}
	\frac{d}{dt} \bigg\{ \beta\io \ueps\ln\ueps 
	+ \frac{\gamma\chi}{2\alpha} \io \frac{|\nabla\veps|^2}{\veps}
	+ \frac{\chi}{2} \io \frac{|\nabla\weps|^2}{\weps} \bigg\} \nn\\
	&\le& - \beta D \io \frac{|\nabla\ueps|^2}{\ueps}
	+ \beta\chi \io \Feps'(\ueps) \nabla\ueps\cdot\nabla\weps
	+ \beta\io \Feps(\ueps) \ln \ueps \cdot \weps
	+ \beta \io \Feps(\ueps) \weps \nn\\
	& & + \frac{\gamma\chi}{2} \io \frac{\weps}{\veps} |\nabla \veps|^2
	+ \gamma\chi \io \nabla\veps\cdot\nabla\weps \nn\\
	& & - \frac{\chi}{2} \io \weps |D^2\ln\weps|^2 \nn\\
	& & - \frac{\beta\chi}{2} \io \frac{\Feps(\ueps)}{\weps} |\nabla\weps|^2
	- \beta\chi \io \Feps'(\ueps) \nabla\ueps\cdot\nabla\weps \nn\\
	& & - \frac{\gamma\chi}{2} \io \frac{\veps}{\weps} |\nabla\weps|^2
	- \gamma\chi \io \nabla\veps\cdot\nabla\weps \nn\\[2mm]
	&\le& - \beta D \io \frac{|\nabla\ueps|^2}{\ueps}
	- \frac{\chi}{2c_1} \io \frac{|\Delta\weps|^2}{\weps}
	- \frac{\chi}{2c_2} \io \frac{|\nabla\weps|^4}{\weps^3} \nn\\
	& & + \frac{\gamma\chi}{2} \io \frac{\weps}{\veps} |\nabla\veps|^2
	+ \beta \io \Feps(\ueps)\ln\ueps \cdot\weps 
	+ \beta\io \Feps(\ueps)\weps
	\qquad \mbox{for all } t>0,
  \eea
  where using that 	%with $\kappa>0$ as in (\ref{kappa}) we have
  \be{4.9}
	\|\weps(\cdot,t)\|_{L^\infty(\Omega)} \le K e^{-\kappa t} \le K
	\qquad \mbox{for all } t>0
  \ee
  by (\ref{winfty}) and (\ref{4.1}), we can estimate
  \be{4.10}
	\frac{\chi}{2c_1} \io \frac{|\Delta\weps|^2}{\weps}
	\ge \frac{\chi}{2c_1 K} \io |\Delta\weps|^2
	\qquad \mbox{and} \qquad
	\frac{\chi}{2c_2} \io \frac{|\nabla\weps|^4}{\weps^3}
	\ge \frac{\chi}{2c_2 K^3} \io |\nabla\weps|^4
  \ee
  for all $t>0$.
  As from (\ref{4.9}) and the validity of the inequality $z\ln z\ge -\frac{1}{e}$ for all $z>0$ we moreover obtain that
  \bas
	\frac{\gamma\chi}{2}
	\io \frac{\weps}{\veps} |\nabla\veps|^2
	&\le& \frac{\gamma\chi K}{2} e^{-\kappa t} \cdot \io \frac{|\nabla \veps|^2}{\veps} \\
	&=& \alpha K e^{-\kappa t} \cdot \bigg\{ \F(t) - \beta \io \ueps\ln\ueps
	- \frac{\chi}{2} \io \frac{|\nabla\weps|^2}{\weps} \bigg\} \\
	&\le& \alpha K e^{-\kappa t} \cdot \F(t)
	- \alpha\beta K e^{-\kappa t} \cdot \int_{\{\ueps<1\}} \ueps\ln\ueps \\
	&\le& \alpha K e^{-\kappa t} \cdot \F(t)
	+ \frac{\alpha\beta K |\Omega|}{e} \cdot e^{-\kappa t}
	\qquad \mbox{for all } t>0
  \eas
  and, similarly,
  \bas
	\beta \io \Feps(\ueps) \ln\ueps \cdot \weps
	&\le& \beta \int_{\{\ueps>1\}} \Feps(\ueps) \ln\ueps \cdot \weps \\
	&\le& \beta K e^{-\kappa t} \cdot \int_{\{\ueps>1\}} \ueps\ln\ueps \\
	&=& \beta K e^{-\kappa t} \cdot \bigg\{ \io \ueps\ln\ueps - \int_{\{\ueps<1\}} \ueps\ln\ueps \bigg\} \\
	&\le& \beta K e^{-\kappa t} \cdot \bigg\{ \io \ueps\ln \ueps + \frac{|\Omega|}{e}\bigg\} \\
	&\le& \beta K e^{-\kappa t} \cdot \F(t)
	+ \frac{\beta K |\Omega|}{e} e^{-\kappa t}
	\qquad \mbox{for all } t>0
  \eas
  as well as
  \bas
	\beta\io \Feps(\ueps)\weps
	\le \beta K e^{-\kappa t} \io \ueps
	\le \beta K \cdot \bigg\{ \io u_0 + \frac{1}{\beta} \io w_0\bigg\}
	\qquad \mbox{for all } t>0
  \eas
  due to (\ref{mass}), we readily infer that (\ref{4.8}) entails (\ref{4.4}) with some appropriately large $C(K)>0$.
\qed
\subsection{Resulting a priori estimates}
Since the factor $e^{-\kappa t}$ appearing on the right of (\ref{4.4}) is integrable over $(0,\infty)$,
our first conclusion from Lemma \ref{lem4} asserts some space-time bounds for the quantities making up
the dissipation rate $\D$ therein even over this entire unbounded time interval.
\begin{lem}\label{lem5}
  Assume that $u_0$ and $v_0$ comply with (\ref{iuv}).
  Then for all $M>0$ and $K>0$ 
  there exists $C(M,K)>0$ such that if $w_0$ satisfies (\ref{iw}) as well as (\ref{M}) and (\ref{4.1}),
  we have
  \be{5.1}
	\int_0^\infty \io \frac{|\nabla\ueps|^2}{\ueps} \le C(M,K)
	\qquad \mbox{for all } \eps\in (0,1)
  \ee
  and
  \be{5.2}
	\int_0^\infty \io |\Delta\weps|^2 \le C(M,K)
	\qquad \mbox{for all } \eps\in (0,1)
  \ee
  as well as
  \be{5.3}
	\int_0^\infty \io |\nabla\weps|^4 \le C(M,K)
	\qquad \mbox{for all } \eps\in (0,1).
  \ee
\end{lem}
\proof
  On invoking Lemma \ref{lem4}, we obtain
  $c_1=c_1(K)>0$ and $c_2=c_2(K)>0$ such that whenever $\eps\in (0,1)$, with $\kappa>0$ as in (\ref{kappa}) we have
  \be{5.4}
	\frac{d}{dt} \F(t) +  + c_1 \D(t) \le c_2 e^{-\kappa t} \cdot (1+\F(t))
	\qquad \mbox{for all } t>0,
  \ee
  from which on dropping the second summand on the left we infer by comparison that
  \bas
	\F(t)
	&\le& \F(0) e^{c_2 \int_0^t e^{-\kappa s} ds}
	+ \int_0^t e^{c_2\int_s^t e^{-\kappa\sigma} d\sigma} \cdot c_2 e^{-\kappa s} ds \\
	&\le& c_3 \equiv c_3(M,K):= \bigg\{ \beta \cdot \bigg| \io u_0\ln u_0 \bigg| 
	+ \frac{\gamma\chi}{2\alpha} \io \frac{|\nabla v_0|^2}{v_0}
	+ \frac{\chi M}{2} \bigg\} \cdot e^\frac{c_2}{\kappa}
	+ \frac{c_2}{\kappa} \cdot e^\frac{c_2}{\kappa}
	\qquad \mbox{for all } t>0
  \eas
  according to (\ref{M}).
  Thereafter, an integration of (\ref{5.4}) shows that
  \bas
	c_1\int_0^t \D(s) ds
	&\le& \F(0) - \F(t)
	+ c_2\int_0^t e^{-\kappa s} \cdot (1+\F(s)) ds \\
	&\le& c_3 -\beta \io \ueps(\cdot,t)\ln\ueps(\cdot,t) 
	+ \frac{c_2\cdot (1+c_3)}{\kappa} \\
	&\le& c_3
	+ \frac{\beta|\Omega|}{e}
	+ \frac{c_2\cdot (1+c_3)}{\kappa} :=C(M,K)
  \eas
  and thereby implies (\ref{5.1})-(\ref{5.3}) in view of the definition (\ref{4.3}) of $\D$.
\qed
When restricted to bounded time intervals and combined with (\ref{mass}), 
by means of straightforward interpolation the above inequalities can be seen to entail 
some estimates also for $\ueps$ and $\nabla\ueps$ themselves, without appearance of any weight function.
\begin{lem}\label{lem6}
  Assume (\ref{iuv}), and let $M>0$, $K>0$ and $T>0$. Then there exists $C(M,K,T)>0$ with the property that
  whenever $w_0$ satisfies (\ref{iw}), (\ref{M}) and (\ref{4.1}),
  we have
  \be{6.1}
	\int_0^T \io \ueps^\frac{n+2}{n} \le C(M,K,T)
	\qquad \mbox{for all } \eps\in (0,1)
  \ee
  and
  \be{6.2}
	\int_0^T \io |\nabla\ueps|^\frac{n+2}{n+1} \le C(M,K,T)
	\qquad \mbox{for all } \eps\in (0,1).
  \ee
\end{lem}
\proof
  The estimate (\ref{6.1}) follows by straightforward interpolation between the spatio-temporal $L^2$ estimate
  for $\nabla \sqrt{\ueps}$ from (\ref{5.1}) and (\ref{mass}). Thereafter, (\ref{6.2}) becomes a consequence
  of (\ref{6.1}) and again (\ref{5.1}) by means of the H\"older inequality (see also \cite[Lemma 3.2]{win_ct_abs}).
\qed
For the second component, on the basis of the ODE therefor in (\ref{0eps}), as well as (\ref{5.3}), 
we obtain the following.
\begin{lem}\label{lem8}
  If (\ref{iuv}) holds, then for all $M>0$, $K>0$ and $T>0$ one can find $C(M,K,T)>0$ such that
  if $w_0$ satisfies (\ref{iw}), (\ref{M}) and (\ref{4.1}), then for all $\eps\in (0,1)$,
  \be{8.1}
	\io |\nabla\veps(\cdot,t)|^4 \le C(M,K,T)
	\qquad \mbox{for all } t\in (0,T).
  \ee
  Moreover,
  \be{8.2}
	|v_{\eps t}(x,t)| \le \alpha K e^\frac{\alpha K}{\kappa} \|v_0\|_{L^\infty(\Omega)}
	\qquad \mbox{for all $x\in\Omega, t>0$ and } \eps\in (0,1),
  \ee
  where $\kappa>0$ is as in (\ref{kappa}).
\end{lem}
\proof
  We differentiate the identity
  \bas
	\veps(x,t) = v_0(x) e^{-\int_0^t \alpha \weps(x,s) ds},
	\qquad x\in\Omega, \ t>0,
  \eas
  and use (\ref{winfty}) and (\ref{4.1}) to estimate
  \bas
	|\nabla\veps(x,t)|
	&=& \bigg| \alpha v_0(x) e^{\alpha \int_0^t \weps(x,s) ds} \cdot \int_0^t \nabla\weps(x,s) ds 
	+ \nabla v_0(x) \cdot e^{\alpha\int_0^t \weps(x,s) ds} \bigg| \\
	&\le& \alpha \|v_0\|_{L^\infty(\Omega)}
	\cdot e^{\alpha K \int_0^t e^{-\kappa s} ds} \int_0^t |\nabla \weps(x,s)| ds 
	+ \|\nabla v_0\|_{L^\infty(\Omega)} \cdot e^{\alpha K \int_0^t e^{-\kappa s} ds} \\
	&\le& c_1 \int_0^t |\nabla \weps(x,s)|ds + c_2
	\qquad \mbox{for all $x\in\Omega$ and } t>0,
  \eas
  where $c_1:=\alpha\|v_0\|_{L^\infty(\Omega)} \cdot e^\frac{\alpha K}{\kappa}$ and 
  $c_2:=\|\nabla v_0\|_{L^\infty(\Omega)} \cdot e^\frac{\alpha K}{\kappa}$.
  After an integration over $\Omega$, according to the H\"older inequality this shows that
  \bas
	\io |\nabla \veps(x,t)|^4 dx
	&\le& 8c_1^4 \io \bigg\{ \int_0^t |\nabla \weps(x,s)| ds \bigg\}^4 dx + 8c_2^4 \\
	&\le& 8c_1^4 \cdot \io \bigg\{ \int_0^t |\nabla\weps(x,s)|^4 ds \bigg\} \cdot t^3 dx + 8c_2^4 \\
	&\le& 8c_1^4 T^3 \cdot \int_0^\infty \io |\nabla\weps|^4 + 8c_2^4
	\qquad \mbox{for all } t>0
  \eas
  and hence implies (\ref{8.1}) due to (\ref{5.3}).\\
  The explicit inequality in (\ref{8.2}) is an immediate consequence of (\ref{winfty}), (\ref{4.1}) and (\ref{vinf}).
\qed
Our derivation of bounds for the time derivatives of $\ueps$ and $\weps$ by means of Lemma \ref{lem5} in the next two lemmata
is rather straightforward.
\begin{lem}\label{lem21}
  Assume (\ref{iuv}). Then for all $M>0$, $K>0$ and $T>0$ there exists $C(M,K,T)>0$ with the property that
  whenever $w_0$ satisfies (\ref{iw}), (\ref{M}) and (\ref{4.1}), we have
  \be{21.1}
	\int_0^T \|u_{\eps t}(\cdot,t)\|_{(W^{1,p}(\Omega))^\star}^q dt \le C(M,K,T)
	\qquad \mbox{for all } \eps\in (0,1),
  \ee
  where
  \be{pq}
	p:=\max\Big\{ n+2 \, , \, \frac{4(n+2)}{6-n} \Big\}	
	\qquad \mbox{and} \qquad
	q:=\min\Big\{ \frac{n+2}{n+1} \, , \, \frac{4(n+2)}{5n+2} \Big\} \, >1.
  \ee
\end{lem}
\proof
  Given $t>0$ and $\psi\in C^1(\bom)$, using (\ref{0eps}) along with (\ref{F1}) and the H\"older inequality we obtain
  \bas
	\bigg| \io u_{\eps t}(\cdot,t)\cdot\psi \bigg|
	&=& \bigg| -D\io \nabla\ueps\cdot\nabla\psi
	+ \chi \io \ueps\Feps'(\ueps) \nabla\weps\cdot\nabla\psi
	+ \io \Feps(\ueps)\weps\psi \bigg| \\
	&\le& D\|\nabla\ueps\|_{L^\frac{n+2}{n+1}(\Omega)} \|\nabla\psi\|_{L^{n+2}(\Omega)}
	+ \chi \|\ueps\|_{L^\frac{n+2}{n}(\Omega)} \|\nabla\weps\|_{L^4(\Omega)} 
		\|\nabla\psi\|_{L^\frac{4(n+2)}{6-n}(\Omega)} \\
	& & + \|\ueps\|_{L^\frac{n+2}{n}(\Omega)} \|\weps\|_{L^\infty(\Omega)} \|\psi\|_{L^\frac{n+2}{2}(\Omega)}
	\qquad \mbox{for all } \eps\in (0,1).
  \eas
  Since (\ref{pq}) warrants that $W^{1,p}(\Omega)$ is continuously embedded into $W^{1,n+2}(\Omega)$,
  $W^{1,\frac{4(n+2)}{6-n}}(\Omega)$ and $L^\frac{n+2}{2}(\Omega)$, we thus infer the existence of $c_1>0$ such that
  \bas
	\|u_{\eps t}(\cdot,t)\|_{(W^{1,p}(\Omega))^\star}
	&\le& c_1 \|\nabla\ueps\|_{L^\frac{n+2}{n+1}(\Omega)}
	+c_1 \|\ueps\|_{L^\frac{n+2}{n}(\Omega)} \|\nabla\weps\|_{L^4(\Omega)}  \\
	& & + c_1\|\ueps\|_{L^\frac{n+2}{n}(\Omega)} \|\weps\|_{L^\infty(\Omega)} 
	\qquad \mbox{for all $t>0$ and } \eps\in (0,1),
  \eas
  so that in view of Young's inequality, an integration yields
  \bas
	\int_0^T \|u_{\eps t}(\cdot,t)\|_{(W^{1,p}(\Omega))^\star}^q dt
	&\le& c_1^q \int_0^T \|\nabla\ueps(\cdot,t)\|_{L^\frac{n+2}{n+1}(\Omega)}^q dt 
	+ c_1^q \int_0^T \|\ueps(\cdot,t)\|_{L^\frac{n+2}{n}(\Omega)}^q \|\nabla\weps(\cdot,t)\|_{L^4(\Omega)}^q dt \\
	& & + c_1^q \int_0^T \|\ueps(\cdot,t)\|_{L^\frac{n+2}{n}(\Omega)}^q \|\weps(\cdot,t)\|_{L^\infty(\Omega)}^q dt \\
	&\le& c_1^q \int_0^T \|\nabla\ueps(\cdot,t)\|_{L^\frac{n+2}{n+1}(\Omega)}^\frac{n+2}{n+1} dt + c_1 T \\
	& & + c_1^q \int_0^T \|\ueps(\cdot,t)\|_{L^\frac{n+2}{n}(\Omega)}^\frac{n+2}{n} dt
	+ c_1^q \int_0^T \|\nabla\weps(\cdot,t)\|_{L^4(\Omega)}^4 dt + c_1 T \\
	& & + c_1 \int_0^T \|\ueps(\cdot,t)\|_{L^\frac{n+2}{n}(\Omega)}^\frac{n+2}{n} dt
	+ c_1^q \int_0^T \|\weps(\cdot,t)\|_{L^\infty(\Omega)}^\frac{n+2}{n+2-nq} dt
%	\qquad \mbox{for all } \eps\in (0,1),
  \eas
  for all $\eps\in (0,1)$,
  because $q\le \frac{n+2}{n+1}<\frac{n+2}{n}$ and $\frac{n}{n+2} + \frac{1}{4} + \frac{1}{q}\le 1$ by (\ref{pq}).
  Therefore, (\ref{21.1}) results from (\ref{6.2}), (\ref{6.1}), (\ref{5.3}) and (\ref{winfty}).
\qed
\begin{lem}\label{lem22}
  Assume (\ref{iuv}). Then for all $M>0$, $K>0$ and $T>0$ there exists $C(M,K,T)>0$ with the property that
  whenever $w_0$ satisfies (\ref{iw}), (\ref{M}) and (\ref{4.1}), we have
  \be{22.1}
	\int_0^T \|w_{\eps t}(\cdot,t)\|_{(W_0^{1,\infty}(\Omega))^\star} dt \le C(M,K,T)
	\qquad \mbox{for all } \eps\in (0,1).
  \ee
\end{lem}
\proof
  We fix $\psi\in C_0^\infty(\Omega)$ such that $\|\psi\|_{W^{1,\infty}(\Omega)} \le 1$ and then see using (\ref{0eps}),
  (\ref{F1}) and Young's inequality that
  \bas
	\bigg| \io w_{\eps t} \cdot\psi\bigg|
	&=& \bigg| - \io \nabla\weps\cdot\nabla\psi
	- \beta\io \Feps(\ueps)\weps\psi
	- \gamma\io \veps\weps\psi \bigg| \\
	&\le& \io |\nabla\weps|
	+ \beta \io \ueps\weps
	+ \gamma\io \veps\weps \\
	&\le& \io |\nabla\weps|^4 + |\Omega|
	+ \beta\|\ueps\|_{L^1(\Omega)} \|\weps\|_{L^\infty(\Omega)}
	+ \gamma \|\veps\|_{L^\infty(\Omega)} \|\weps\|_{L^\infty(\Omega)} |\Omega|
  \eas
  and that therefore
  \bas
	\int_0^T \|w_{\eps t}(\cdot,t)\|_{(W_0^{1,\infty}(\Omega))^\star} dt
	&\le& \int_0^T \io |\nabla\weps|^4 + |\Omega| T \\
	& & + \beta \cdot \bigg\{ \io u_0 + \frac{1}{\beta} \io w_0\bigg\} \cdot \|w_0\|_{L^\infty(\Omega)} \cdot T \\
	& & + \gamma \|v_0\|_{L^\infty(\Omega)} e^{\frac{\alpha}{\kappa} \|w_0\|_{L^\infty(\Omega)}}
	\cdot \|w_0\|_{L^\infty(\Omega)} \cdot |\Omega| T
	\qquad \mbox{for all } \eps\in (0,1)
  \eas
  thanks to (\ref{mass}), (\ref{winfty}) and (\ref{vinf}).
  According to (\ref{5.3}) and (\ref{4.1}), this hence entails (\ref{22.1}).
\qed
Our collection of estimates now enables us to pass to the limit in the following sense:
%We have now collected all ingredients sufficient for passing to the limit in the following sense:
%
%
\begin{lem}\label{lem9}
  Let $n\le 5$, and suppose that $u_0, v_0$ and $w_0$ satisfy (\ref{iuv}) and (\ref{iw}).
  Then there exist $(\eps_j)_{j\in\N}\subset (0,1)$ and nonnegative functions
  \be{9.1}
	\left\{ \begin{array}{l}
	u \in L^\infty((0,\infty);L^1(\Omega)) 
	\cap L^\frac{n+2}{n}_{loc}([0,\infty);L^\frac{n+2}{n}(\Omega))
	\cap L^\frac{n+2}{n+1}_{loc}([0,\infty);W^{1,\frac{n+2}{n+1}}(\Omega)), \\[1mm]
	v\in L^\infty(\Omega\times (0,\infty))
	\qquad \mbox{and} \\[1mm]
	w\in L^\infty(\Omega\times (0,\infty)) \cap L^4_{loc}([0,\infty);W^{1,4}(\Omega))
	\cap L^2_{loc}([0,\infty);W^{2,2}(\Omega))
	\end{array} \right.
  \ee
  such that $\eps_j\searrow 0$ as $j\to\infty$ and that
  \bas
	& & \ueps\to u
	\qquad \mbox{a.e.~in } \Omega\times (0,\infty), 
	\label{9.2} \\
	& & \ueps\to u, \Feps(\ueps)\to u
	\quad \mbox{and} \quad \ueps \Feps'(\ueps) \to u
	\qquad \mbox{in } L^p_{loc}(\bom\times [0,\infty))
	\quad \mbox{for all } p\in \Big[1,\frac{n+2}{n}\Big),
	\label{9.3} \\
	& & \nabla\ueps \wto \nabla u
	\qquad \mbox{in } L^\frac{n+2}{n+1}_{loc}(\bom\times [0,\infty)),
	\label{9.4} \\
	& & u_{\eps t} \wto u_t
	\qquad \mbox{in } L^q_{loc}([0,\infty);(W^{1,p}(\Omega))^\star)
	\quad \mbox{for $p$ and $q$ as in (\ref{pq}),}
	\label{9.44} \\
	& & \veps\to v
	\quad \mbox{a.e.~in } \Omega\times (0,\infty)
	\qquad \mbox{and} \qquad
	\veps(\cdot,t) \to v(\cdot,t)
	\quad \mbox{a.e.~in $\Omega$ for a.e. } t>0,
	\label{9.5} \\
	& & \weps\to w
	\qquad \mbox{a.e.~in } \Omega\times (0,\infty), 
	\label{9.6} \\
	& & \nabla\weps \wto \nabla w
	\qquad \mbox{in } L^4_{loc}(\bom\times [0,\infty))
	\qquad \mbox{and}
	\label{9.7} \\
	& & \Delta \weps \wto \Delta w
	\qquad \mbox{in } L^2_{loc}(\bom\times [0,\infty))
	\label{9.8}
  \eas
  as $\eps=\eps_j\searrow 0$. Moreover, $(u,v,w)$ is a global weak solution of (\ref{0}) in the sense of Definition
  \ref{defi_weak}.
\end{lem}
\proof
  In view of the estimates obtained in Lemmata \ref{lem5}, \ref{lem6}, \ref{lem8}, \ref{lem21} and \ref{lem22},
  this can be achieved by means of essentially straightforward extraction procedures based on the Aubin-Lions lemma
  and the Vitali convergence theorem; a corresponding reasoning in a closely related situation can be found in
  \cite{win_ct_abs}.
\qed
As a by-product, we immediately obtain our main result on global solvability in the original problem (\ref{0}).\abs
\proofc of Theorem \ref{theo222}. \quad
  This directly results from Lemma \ref{lem9}.
\qed
\mysection{Migration-driven benefit}\label{sect3}
%\mysection{Qualitative analysis: migration-driven benefit}
%
%
%
%
%
%
%
%
In order to derive the qualitative properties claimed in Theorem \ref{theo18}, let us first draw some essentially evident
consequences of our estimates gained above.
\begin{lem}\label{lem10}
  Assume (\ref{iuv}). Then for all $M>0$, $K>0$ and $T>0$ there exists $C(M,K,T)>0$ with the property that
  whenever $w_0$ satisfies (\ref{iw}), (\ref{M}) and (\ref{4.1}),
  the global weak solution $(u,v,w)$ of (\ref{0}) from Lemma \ref{lem9} has the properties that
  \be{10.1}
	\int_0^T \io u^\frac{n+2}{n} \le C(M,K,T)
  \ee
  and
  \be{10.2}
	\int_0^T \io |\nabla u|^\frac{n+2}{n+1} \le C(M,K,T)
  \ee
  as well as
  \be{10.3}
	\int_0^T \io |\nabla w|^4 \le C(M,K,T)
  \ee
  and
  \be{10.4}
	\int_0^T \|u_t(\cdot,t)\|_{(W^{1,p}(\Omega))^\star}^q dt \le C(M,K,T)
  \ee
  with $p>1$ and $q>1$ taken from (\ref{pq}).
\end{lem}
\proof
  In view of (\ref{9.2}) and Fatou's lemma, the estimate in (\ref{10.1}) directly results from (\ref{6.1}).
  Likewise, (\ref{10.2}), (\ref{10.3}) and (\ref{10.4}) are consequences of (\ref{6.2}), (\ref{5.3}) and (\ref{21.1})
  when combined with (\ref{9.4}), (\ref{9.7}) and (\ref{9.44}) by means of a standard argument based on lower
  semicontinuity with respect to weak convergence.
\qed
\subsection{Detecting a genuine Lyapunov functional}
Now a key toward our derivation of Theorem \ref{theo18} will consist in the detection of another Lyapunov-type preoperty
of (\ref{0eps}) which, unlike that from Lemma \ref{lem4}, will become manifest in a genuine energy inequality.
This will be derived in Lemma \ref{lem12} on adequately combining the following simple observations.
\begin{lem}\label{lem11}
  Suppose that (\ref{iuv}) and (\ref{iw}) hold. Then for all $\eps\in (0,1)$,
  \be{11.1}
	- \frac{d}{dt} \io \ln \ueps
	+ D \io \frac{|\nabla\ueps|^2}{\ueps^2}
	= \chi \io \Feps'(\ueps) \nabla\ueps\cdot\nabla\weps
	- \io \frac{\Feps(\ueps)}{\ueps} \cdot \weps
	\qquad \mbox{for all } t>0
  \ee
  and
  \be{11.11}
	\frac{1}{2} \frac{d}{dt} \io \weps^2 + \io |\nabla\weps|^2
	= - \beta \io \Feps(\ueps) \weps^2 - \gamma \io \veps\weps^2
	\qquad \mbox{for all } t>0
  \ee
  as well as
  \be{11.2}
	\int_0^t \io |\nabla\weps|^2 \le \frac{1}{2} \io w_0^2
	\qquad \mbox{for all } t>0.
  \ee
\end{lem}
\proof
  The identities in (\ref{11.1}) and (\ref{11.11}) immediately result from testing the first and the third equation
  in (\ref{0eps}) by $\frac{1}{\ueps}$ and $\weps$, respectively.
  Thereafter, (\ref{11.2}) follows from (\ref{11.11}) on integrating in time and dropping two favorably signed summands.
\qed
Taking suitable linear combinations of the above identities and additionally using further basic
information on mass evolution in (\ref{0eps}), we indeed obtain the following crucial inequality.
\begin{lem}\label{lem12}
  Assume (\ref{iuv}) and (\ref{iw}), and with $\kappa>0$ as in (\ref{kappa}), let
  \be{ab}
	a:=\frac{\alpha+\frac{1}{4}}{\kappa}
	\qquad \mbox{and} \qquad
	b:=\frac{\chi^2}{4D}.
  \ee
  Then for all $\eps\in (0,1)$,
  \be{12.2}
	\frac{d}{dt} \bigg\{ \io \ln \frac{\veps}{\ueps} + a\io \weps + b\io \weps^2 \bigg\}
	\le - \frac{D}{2}\io \frac{|\nabla\ueps|^2}{\ueps^2} - \frac{1}{4} \io \weps
	\qquad \mbox{for all } t>0.
  \ee
\end{lem}
\proof
  Using (\ref{0eps}), for $t>0$ we compute
  \bas
	\frac{d}{dt} \io \ln \veps = \alpha\io \weps
	\qquad \mbox{for all } t>0
  \eas
  and 
  \bas
	\frac{d}{dt} \io \weps
	= - \beta \io \Feps(\ueps)\weps - \gamma \io \veps\weps
	\qquad \mbox{for all } t>0,
  \eas
  whence in view of (\ref{vinf}) and (\ref{kappa}) we have
  \bas
	a\frac{d}{dt} \io \weps 
	\le - a\kappa \io \weps
	\qquad \mbox{for all } t>0.
  \eas
  As moreover
  \bas
	b\frac{d}{dt} \io \weps^2 \le -2b \io |\nabla\weps|^2
	\qquad \mbox{for all } t>0
  \eas
  by (\ref{11.11}), recalling (\ref{11.1}) we obtain
  \bea{12.3}
%	\hspace*{-20mm}
	\frac{d}{dt} \bigg\{ \io \ln \frac{\veps}{\ueps} + a\io \weps + b\io \weps^2 \bigg\} 	%\nn\\
	&\le& \alpha\io \weps \nn\\
	& & - D\io \frac{|\nabla\ueps|^2}{\ueps^2}
	+ \chi \io \frac{\Feps'(\ueps)}{\ueps} \nabla\ueps\cdot\nabla\weps 
	- \io \frac{\Feps(\ueps)}{\ueps} \weps \nn\\
	& & - a\kappa \io \weps \nn\\
	& & -2b\io |\nabla\weps|^2 \nn\\
	&\le& -D \io \frac{|\nabla\ueps|^2}{\ueps^2}
	+ \chi \io \frac{\Feps'(\ueps)}{\ueps} \nabla\ueps\cdot\nabla\weps \nn\\
	& & -2b\io |\nabla\weps|^2 - \frac{1}{4} \io \weps
	\qquad \mbox{for all } t>0,
  \eea
  because $\Feps\ge 0$ by (\ref{F1}), and because $\alpha-a\kappa=-\frac{1}{4}$ due to (\ref{ab}).
  Here we use Young's inequality and again (\ref{F1}) to estimate
  \bas
	\chi \io \frac{\Feps'(\ueps)}{\ueps} \nabla\ueps\cdot\nabla\weps 
	&\le& \frac{D}{2} \io \frac{|\nabla\ueps|^2}{\ueps^2}
	+ \frac{\chi^2}{2D} \io \Feps'^2(\ueps) |\nabla\weps|^2 \\
	&\le& \frac{D}{2} \io \frac{|\nabla\ueps|^2}{\ueps^2}
	+ \frac{\chi^2}{2D} \io |\nabla\weps|^2
	\qquad \mbox{for all } t>0,
  \eas
  so that (\ref{12.2}) becomes a consequence of (\ref{12.3}) and the fact that $\frac{\chi^2}{2D}=2b$ by (\ref{ab}).
\qed
By integration, this immediately implies the following.
\begin{lem}\label{lem122}
  Assume (\ref{iuv}) and (\ref{iw}), and let $a>0$ and $b>0$ be as in (\ref{ab}). Then for all $\eps\in (0,1)$,
  \be{122.1}
	\io \ln \frac{\veps(\cdot,t)}{\ueps(\cdot,t)}
	+ \frac{D}{2} \int_0^t \io \frac{|\nabla\ueps|^2}{\ueps^2} + \frac{1}{4} \int_0^t\io \weps
	\le \io \ln \frac{v_0}{u_0} + a\io w_0 + b \io w_0^2
	\qquad \mbox{for all } t>0.
  \ee
  In particular, whenever (\ref{iuv}) holds and $M>0$ and $K>0$ are given, one can find 
  $C(M,K)>0$ such that if $w_0$ satisfies (\ref{iw}), (\ref{M}) and (\ref{4.1}), then
  for all $\eps\in (0,1)$,
  \be{122.2}
	\io \ln \ueps(\cdot,t) \ge -C(M,K)
	\qquad \mbox{for all } t>0
  \ee
  and
  \be{122.3}
	\int_0^\infty \io \frac{|\nabla\ueps|^2}{\ueps^2} \le C(M,K)
  \ee
\end{lem}
\proof
  According to Lemma \ref{lem12}, we have
  \bas
	\frac{d}{dt} \bigg\{ \io \ln \frac{\veps}{\ueps} + a\io \weps + b\io \weps^2 \bigg\}
	\le - \frac{D}{4} \io \frac{|\nabla\ueps|^2}{\ueps^2} 
	\qquad \mbox{for all } t>0,
  \eas
  from which (\ref{122.1}) and hence also (\ref{122.2}) and (\ref{122.3}) immediately result upon integration.
\qed	
\subsection{Estimating migration effects from below. Proof of Theorem \ref{theo18}}
It will now be of fundamental importance to make sure that the dissipative action expressed through the appearance of the 
second summand on the left of (\ref{122.1}) is conveniently substantial. 
At the level of approximate solutions, a preparation for this can be gained by once again going back to (\ref{11.1}):
\begin{lem}\label{lem13}
  Assume (\ref{iuv}) and (\ref{iw}). Then for all $\eps\in (0,1)$,
  \bea{13.1}
	2D \int_0^{t_\star} \io \frac{|\nabla\ueps|^2}{\ueps^2}
	\ge \io \ln\ueps(\cdot,t_\star) 
	- \io \ln u_0
	- \frac{\chi^2}{8D} \io w_0^2
	- \int_0^{t_\star} \io \weps
	\qquad \mbox{for all } t_\star>0.
  \eea
\end{lem}
\proof
  Again starting from (\ref{11.1}), we now estimate the first summand on the right-hand side therein from below
  by using Young's inequality and (\ref{F1}) according to
  \bas
	\chi \io \frac{\Feps'(\ueps)}{\ueps} \nabla\ueps\cdot\nabla\weps
	&\ge& -D \io \frac{|\nabla\ueps|^2}{\ueps^2}
	- \frac{\chi^2}{4D} \io \Feps'^2(\ueps) |\nabla\weps|^2 \\
	&\ge& -D \io \frac{|\nabla\ueps|^2}{\ueps^2}
	- \frac{\chi^2}{4D} \io |\nabla\weps|^2
	\qquad \mbox{for all } t>0,
  \eas
  whence integrating (\ref{11.1}) in time we obtain
  \bas
	2D\int_0^{t_\star} \io \frac{|\nabla\ueps|^2}{\ueps^2}
	\ge \io \ln\ueps(\cdot,t_\star)
	- \io \ln u_0
	- \frac{\chi^2}{4D} \int_0^{t_\star} \io |\nabla\weps|^2
	- \int_0^{t_\star} \io \frac{\Feps(\ueps)}{\ueps} \weps
	\quad \mbox{for all } t_\star>0.
  \eas
  As $\frac{\Feps(\ueps)}{\ueps} \le 1$ by (\ref{F1}), in view of (\ref{11.2}) this entails (\ref{13.1}).
\qed
In order to appropriately pass to the limit $\eps\searrow 0$ in the latter inequalities, 
let us make sure that also the potentially
singular quantity $\ln \ueps$ appearing therein exhibits a favorably controlable behavior when $\eps$ becomes small.
This will once more be achieved by means of an Aubin-Lions type argument on the basis of the following additional
regularity property.
\begin{lem}\label{lem123}
  If (\ref{iuv}) holds, then for all $M>0$, $K>0$ and $T>0$ one can find $C(M,K,T)>0$ such that
  whenever $w_0$ satisfies (\ref{iw}), (\ref{M}) and (\ref{4.1}), we have
% and (\ref{iw}) hold, then for all $T>0$ one can find $C(T)>0$ such that
  \be{123.1}
	\int_0^T \Big\| \partial_t \ln \ueps(\cdot,t) \Big\|_{(W_0^{1,\infty}(\Omega))^\star} dt \le C(M,K,T)
	\qquad \mbox{for all } \eps\in (0,1).
  \ee
\end{lem}
\proof
  For $\psi\in C_0^\infty(\Omega)$ fulfilling $\|\psi\|_{W^{1,\infty}(\Omega)} \le 1$, from the first equation in
  (\ref{0eps}) we see that for all $t>0$ and $\eps\in (0,1)$,
  \bas
	\bigg| \io \partial_t \ln\ueps(\cdot,t)\cdot\psi \bigg|
	&=& \bigg| \io \frac{1}{\ueps} \cdot \Big\{ D\Delta\ueps
	- \chi \nabla\cdot (\ueps\Feps'(\ueps)\nabla\weps) + \Feps(\ueps) \weps \Big\} \cdot \psi \bigg| \\
	&=& \bigg| D \io \frac{|\nabla\ueps|^2}{\ueps^2} \psi
	- D \io \frac{1}{\ueps} \nabla\ueps\cdot\nabla\psi \\
	& & \quad - \chi \io \frac{\Feps'(\ueps)}{\ueps} (\nabla\ueps\cdot\nabla\weps) \psi
	+ \chi \io \Feps'(\ueps)\nabla\weps\cdot\nabla\psi 
	+ \io \frac{\Feps(\ueps)}{\ueps} \weps\psi \bigg| \\
	&\le& D \io \frac{|\nabla\ueps|^2}{\ueps^2}
	+ D \io \frac{|\nabla\ueps|}{\ueps} 
	+ \chi \io \frac{|\nabla\ueps|}{\ueps} \cdot |\nabla\weps|
	+ \chi \io |\nabla\weps|
	+ \io \weps \\
	&\le& D\io \frac{|\nabla\ueps|^2}{\ueps^2}
	+ \frac{D}{2} \io \frac{|\nabla\ueps|^2}{\ueps^2} + \frac{D|\Omega|}{2} \\
	& & + \frac{\chi}{2} \io \frac{|\nabla\ueps|^2}{\ueps^2}
	+ \frac{\chi}{2} \io |\nabla\weps|^2
	+ \frac{\chi}{2} \io |\nabla\weps|^2 + \frac{\chi|\Omega|}{2} + |\Omega| \cdot \|\weps\|_{L^\infty(\Omega)}
  \eas
  according to (\ref{F1}) and Young's inequality. Thus, for all $t>0$ and $\eps\in (0,1)$ we have
  \bas
	\Big\| \partial_t \ln \ueps(\cdot,t)\Big\|_{(W_0^{1,\infty}(\Omega))^\star}
	\le \frac{3D+\chi}{2} \io \frac{|\nabla\ueps|^2}{\ueps^2}
	+ \chi \io |\nabla\weps|^2
	+ |\Omega|\cdot \|\weps\|_{L^\infty(\Omega)}
	+ \frac{(D+\chi)|\Omega|}{2},
%	\qquad \mbox{for all $t>0$ and } \eps\in (0,1),
  \eas
  which by an integration results in (\ref{123.1}) due to (\ref{122.3}), (\ref{11.2}) and (\ref{winfty}).
\qed
Indeed, we therefore obtain the following.
\begin{lem}\label{lem124}
  i) \ Assume that (\ref{iuv}) and (\ref{iw}) hold, and let $(\eps_j)_{j\in\N}$ be as in Lemma \ref{lem9}.
  Then $\ln u \in L^1(\Omega\times (0,\infty))$, and there exist a null set $N\subset (0,\infty)$ and a subsequence
  $(\eps_{j_l})_{l\in\N}$ of $(\eps_j)_{j\in\N}$ such that
  \be{124.01}
	\ln u(\cdot,t) \in L^1(\Omega)
	\qquad \mbox{for all } t\in (0,\infty)\setminus N
  \ee
  and
  \be{124.1}
	\io \ln \ueps(\cdot,t) \to \io \ln u(\cdot,t)
	\qquad \mbox{for all } t\in (0,\infty)\setminus N
  \ee
  as  $\eps=\eps_{j_l}\searrow 0$.\abs
  ii) \ Assume (\ref{iuv}), and suppose that $(w_{0k})_{k\in\N} \subset W^{1,\infty}(\Omega)$ is such that
  $w_{0k} \ge 0$ in $\Omega$ as well as $\sqrt{w_{0k}} \in W^{1,2}(\Omega)$ for all $k\in\N$ and that
  \be{124.11}
	\sup_{k\in\N} \|w_{0k}\|_{L^\infty(\Omega)} <\infty
	\qquad \mbox{and} \qquad
	\sup_{k\in\N} \io \frac{|\nabla w_{0k}|^2}{w_{0k}} <\infty.
  \ee
  Then the corrseponding global weak solutions $(u_k,v_k,w_k)$ obtained in Theorem \ref{theo222} with
  initial data $(u_0,v_0,w_{0k})$ have the property that $\ln u_k \in L^1_{loc}(\bom\times [0,\infty))$ for all $k\in\N$
  and that $(\ln u_k)_{k\in\N}$ is relatively compact in $L^1_{loc}(\bom\times [0,\infty))$.
\end{lem}
\proof
  i) \ Since  
  \bas
	\io \Big|\ln \ueps(\cdot,t)\Big|
	&=& 2\int_{\{\ueps(\cdot,t)\ge 1\}} \ln \ueps(\cdot,t)
	- \io \ln \ueps(\cdot,t) \\
	&\le& 2\int_{\{\ueps(\cdot,t)\ge 1\}} \ueps(\cdot,t)
	- \io \ln \ueps(\cdot,t) \\
	&\le& 2 \cdot \bigg\{ \io u_0 + \frac{1}{\beta} \io w_0\bigg\}
	- \io \ln \ueps(\cdot,t) 
	\qquad \mbox{for all $t>0$ and } \eps\in (0,1)
  \eas
  by (\ref{mass}), it readily follows from (\ref{122.2}), (\ref{122.3}) and standard $L^1-L^2$ interpolation based on 
  the Gagliardo-Nirenberg inequality that for each $T>0$,
  $(\ln \ueps)_{\eps\in (0,1)} \mbox{is bounded in } L^2((0,T);W^{1,2}(\Omega))$.
%  \bas
%	(\ln \ueps)_{\eps\in (0,1)}
%	\quad \mbox{is bounded in } L^2((0,T);W^{1,2}(\Omega)).
%  \eas
  As Lemma \ref{lem123} asserts that moreover
  $(\partial_t \ln \ueps)_{\eps\in (0,1)} \mbox{is bounded in } L^1((0,T);(W_0^{1,\infty}(\Omega))^\star)$
%  \bas
%	(\partial_t \ln \ueps)_{\eps\in (0,1)}
%	\quad \mbox{is bounded in } L^1((0,T);(W_0^{1,\infty}(\Omega))^\star)
%  \eas
  for any such $T$, an Aubin-Lions lemma warrants strong precompactness of $(\ln \ueps)_{\eps\in (0,1)}$
  in $L^1_{loc}(\bom\times [0,\infty))$, implying that for a suitable subsequence $(\eps_{j_l})_{l\in\N}$
  of $(\eps_j)_{j\in\N}$ we can achieve that
  \be{124.2}
	\ln \ueps \to z
	\qquad \mbox{in }
	L^1_{loc}(\bom\times [0,\infty))
	\mbox{ and a.e.~in } \Omega\times (0,\infty)
  \ee
  and
  \be{124.3}
	\ln \ueps(\cdot,t) \to z(\cdot,t)
	\qquad \mbox{in } L^1(\Omega)
	\quad \mbox{for a.e.~} t>0
  \ee
  as $\eps=\eps_{j_l}\searrow 0$.
  Since (\ref{124.3}) together with (\ref{9.2}) requires that $z=u$ a.e.~in $\Omega\times (0,\infty)$,
  both the inclusion $\ln u\in L^1_{loc}(\bom\times [0,\infty))$ and (\ref{124.1}) become consequences of (\ref{124.2})
  and (\ref{124.3}) if $N\subset (0,\infty)$ is chosen appropriately.\abs
  ii) \ As (\ref{124.11}) warrants that we can achieve independence of the estimates in Lemma \ref{lem122} and
  Lemma \ref{lem123} from $k\in\N$, both properties can by verified by essentially repeating the argument from i).
\qed
When applied to the inequalities from Lemma \ref{lem122} and Lemma \ref{lem13}, this implies an upper bound
for the expression $\io \ln \frac{u}{v}$ involving, inter alia, the difference $\io \ln u - \io \ln u_0$
when evaluated at an arbitrary but fixed time $t_\star$:
\begin{lem}\label{lem14}
  Suppose that (\ref{iuv}) and (\ref{iw}) hold, and let $a>0$ and $b>0$ be as in (\ref{ab}).
  Then there exists a null set $N\subset (0,\infty)$ such that
  \be{14.01}
	\ln \frac{v(\cdot,t)}{u(\cdot,t)} \in L^1(\Omega)
	\qquad \mbox{for all } t\in (0,\infty)\setminus N
  \ee
  and
  \bea{14.1}
	\io \ln \frac{v(\cdot,t)}{u(\cdot,t)}
	&\le& \io \ln \frac{v_0}{u_0}
	+ a\io w_0 + \Big(b+\frac{\chi^2}{32D}\Big) \io w_0^2 
	- \frac{1}{4} \bigg\{ \io \ln u(\cdot,t_\star) - \io \ln u_0 \bigg\} \nn\\[2mm]
	& & \hspace*{40mm}
	\eea
	for all $t_\star\in (0,\infty)\setminus N$ and any $ t\in (t_\star,\infty)\setminus N$.
\end{lem}
\proof
  On combining Lemma \ref{lem122} with Lemma \ref{lem13}, we obtain that whenever $t_\star>0$ and $t>t_\star$,
  \bea{14.2}
	\io \ln \frac{\veps(\cdot,t)}{\ueps(\cdot,t)}
	&\le& \io \ln \frac{v_0}{u_0}
	+ a \io w_0 + b\io w_0^2
	- \frac{D}{2} \int_0^t \io \frac{|\nabla \ueps|^2}{\ueps^2}
	- \frac{1}{4} \int_0^t \io \weps \nn\\
	&\le& \io \ln \frac{v_0}{u_0}
	+ a \io w_0 + b\io w_0^2
	- \frac{D}{2} \int_0^{t_\star} \io \frac{|\nabla \ueps|^2}{\ueps^2}
	- \frac{1}{4} \int_0^{t_\star} \io \weps \nn\\
	&\le& \io \ln \frac{v_0}{u_0}
	+ a \io w_0 + b\io w_0^2 \nn\\
	& & \quad\quad\quad\quad- \frac{1}{4}\left\{\io \ln \ueps(\cdot,t_\star) - \io \ln u_0    \right\} + \frac{\chi^2}{32D} \io w_0^2
  \eea
  for all $\eps\in (0,1)$.
  Now in accordance with Lemma \ref{lem124} we can pick a null set $N_1\subset (0,\infty)$ and a subsequence
  $(\eps_{j_l})_{l\in\N}$ of the sequence $(\eps_j)_{j\in\N}$ from Lemma \ref{lem9} such that
  \be{14.22}
	\ln u(\cdot,t) \in L^1(\Omega)
	\qquad \mbox{for all } t\in (0,\infty)\setminus N_1
  \ee
  and
  \be{14.3}
	\io \ln \ueps(\cdot,t) \to \io \ln u(\cdot,t)
	\quad \mbox{for all } t\in (0,\infty)\setminus N_1
	\qquad \mbox{as } \eps=\eps_{j_l} \searrow 0.
  \ee
  Moreover, from Lemma \ref{lem9} we obtain a null set $N_2 \subset (0,\infty)$ such that as $\eps=\eps_j\searrow 0$
  we have $\veps(\cdot,t)\to v(\cdot,t)$ a.e.~in $\Omega$ for all $t\in (0,\infty)\setminus N_2$.
  Since Lemma \ref{lem2} says that with some $c_1>0$ and $c_2>0$,
  \bas
	c_1 \le \veps(x,t) \le c_2
	\qquad \mbox{for all $(x,t)\in \Omega\times (0,\infty)$ and each } \eps\in (0,1),
  \eas
  due to the dominated convergence theorem this entails that
  \bas
	\ln v(\cdot,t) \in L^1(\Omega)
	\qquad \mbox{for all } t\in (0,\infty)\setminus N_2
  \eas
  and
  \be{14.5}
	\io \ln \veps(\cdot,t) \to \io \ln v(\cdot,t)
	\quad \mbox{for all } t\in (0,\infty)\setminus N_2
	\qquad \mbox{as } \eps=\eps_j \searrow 0
  \ee
  and hence, by (\ref{14.3}), that if we let $N:=N_1 \cup N_2$ and recall (\ref{14.22}), 
  then (\ref{14.01}) holds as well as
  \bas
	\io \ln \frac{\veps(\cdot,t)}{\ueps(\cdot,t)}
	\to \io \ln \frac{v(\cdot,t)}{u(\cdot,t)}
	\quad \mbox{for all } t\in (0,\infty)\setminus N
	\qquad \mbox{as } \eps=\eps_{j_l} \searrow 0.
  \eas
  Once more using (\ref{14.3}), we can thus infer (\ref{14.1}) from (\ref{14.2}).
\qed
Now in order to suitably estimate this difference $\io \ln u(\cdot,t_\star) - \io \ln u_0$ from below at some $t_\star>0$, 
we shall 
invoke a perturbation argument which at its core resorts to a corresponding property of solutions to the heat equation.
\begin{lem}\label{lem15}
  Suppose that $u_0\in C^0(\bom)$ is nonnegative and such that $\io \ln u_0>-\infty$ and $u_0\not\equiv const.$.
  Then there exist $L>0$ and $t_0>0$ such that the solution
  $U\in C^0(\bom\times [0,\infty))\cap C^{2,1}(\bom\times (0,\infty))$ of
  \be{H}
    	\left\{ \begin{array}{rcll}
	U_t &=& D\Delta U,
	\qquad & x\in\Omega, \ t>0, \\[1mm]
	\frac{\partial U}{\partial\nu} &=& 0,
	\qquad & x\in\pO, \ t>0, \\[1mm]
	U(x,0) &=& u_0(x),
	\qquad & x\in\Omega,
	\end{array} \right.
  \ee
  satisfies $U>0$ in $\bom\times [0,\infty)$ and
  \be{15.1}
	\io \ln U(\cdot,t) - \io \ln u_0 \ge L
	\qquad \mbox{for all } t\ge t_0.
  \ee
\end{lem}
\proof
  Since $u_0\not\equiv const.$, a sharp version of Jensen's inequality (see the Appendix below) states that
  $c_1:=\ln \Big\{ \frac{1}{|\Omega|} \io u_0 \Big\} - \frac{1}{|\Omega|} \io \ln u_0$
%  \bas
%	c_1:=\ln \bigg\{ \frac{1}{|\Omega|} \io u_0 \bigg\} - \frac{1}{|\Omega|} \io \ln u_0
%  \eas
  is positive. Using that according to a well-known stabilization property of (\ref{H}) we have
  $U(\cdot,t) \to \frac{1}{|\Omega|} \io u_0$ in $L^\infty(\Omega)$ as $t\to\infty$
%  \bas
%	U(\cdot,t) \to \frac{1}{|\Omega|} \io u_0
%	\quad \mbox{in } L^\infty(\Omega)
%	\qquad \mbox{as } t\to\infty
%  \eas
  and hence, in particular,
  \bas
	\io \ln U(\cdot,t) - \io \ln u_0
	\to \io \ln \bigg\{ \frac{1}{|\Omega|} \io u_0 \bigg\} - \io \ln u_0
	= c_1 |\Omega|
	\qquad \mbox{as } t\to\infty,
  \eas
  we readily conclude that (\ref{15.1}) holds for some suitably large $t_0>0$ if we let, e.g.,
  $L:=\frac{c_1 |\Omega|}{2}$.
\qed
Indeed, the first solution components lie conveniently close to solutions of (\ref{H}) when $w_0$ 
becomes small in an appropriate sense:
\begin{lem}\label{lem16}
  Assume (\ref{iuv}), and suppose that $(w_{0k})_{k\in\N} \subset W^{1,\infty}(\Omega)$ is such that $w_{0k} \ge 0$
  in $\Omega$ and $\sqrt{w_{0k}} \in W^{1,2}(\Omega)$ with
  \be{16.1}
	w_{0k} \to 0
	\quad \mbox{in } L^\infty(\Omega)
	\qquad \mbox{as } k\to\infty
  \ee
  and
  \be{16.2}
	\sup_{k\in\N} \io \frac{|\nabla w_{0k}|^2}{w_{0k}} < \infty.
  \ee
  Moreover, for $k\in\N$ let $(u_k,v_k,w_k)$ denote the global weak solution of (\ref{0}) from Theorem \ref{theo222}
  corresponding to $w_0:=w_{0k}$, and let $U$ solve (\ref{H}).
  Then there exist a null set $N\subset (0,\infty)$ and a subsequence $(u_{k_l})_{l\in\N}$ of
  $(u_k)_{k\in\N}$ such that
  \be{16.3}
	\ln u_{k_l} (\cdot,t) \in L^1(\Omega)
	\qquad \mbox{for all $t\in (0,\infty)\setminus N$ and any } l\in\N,
  \ee
  and that
  \be{16.4}
	\io \ln u_{k_l}(\cdot,t)
	\to \io \ln U(\cdot,t)
	\quad \mbox{for all } t\in (0,\infty)\setminus N
	\qquad \mbox{as } l\to\infty.
  \ee
\end{lem}
\proof
  According to (\ref{16.2}), we may invoke Lemma \ref{lem10} to infer from (\ref{10.1}), (\ref{10.2}) and (\ref{10.4})
  upon another application of the Aubin-Lions lemma that there exist a nonnegative function
  $U\in L^\frac{n+2}{n+1}_{loc}([0,\infty);W^{1,\frac{n+2}{n+1}}(\Omega))$ and a subsequence $(u_{k_l})_{l\in\N}$
  of $(u_k)_{k\in\N}$ such that as $l\to\infty$ we have
  \be{16.5}
	u_{k_l} \to U
	\quad \mbox{in } L^1_{loc}(\bom\times [0,\infty))
	\qquad \mbox{and} \qquad
	\nabla u_{k_l} \wto \nabla u
	\quad \mbox{in } L^1_{loc}(\bom\times [0,\infty)),
  \ee
  where in view of Lemma \ref{lem124} ii), thanks to (\ref{16.1}) and (\ref{16.2}) we can also achieve
  on passing to a further subsequence if necessary that with some null set $N\subset (0,\infty)$ we have
  (\ref{16.3}) as well as
  \be{16.55}
	\io \ln u_{kl}(\cdot,t) \to \io \ln \hatu(\cdot,t)
	\quad \mbox{for all } t\in (0,\infty)\setminus N
	\qquad \mbox{as } l\to\infty.
  \ee
  Since moreover
  \be{16.6}
	w_k 
	\to 0
	\quad \mbox{in } L^\infty(\Omega\times (0,\infty))
	\qquad \mbox{as } k\to\infty
  \ee
  by (\ref{winfty}) and (\ref{16.1}), it follows that if we fix any $\varphi\in C_0^\infty(\bom\times [0,\infty))$
  additionally satisfying $\frac{\partial\varphi}{\partial\nu}=0$ on $\pO\times (0,\infty)$, then in the identity
  \be{16.66}
	- \int_0^\infty u_k \varphi_t
	- \io u_0 \varphi(\cdot,0)
	= -D \int_0^\infty \io \nabla u_k \cdot \nabla \varphi
	+ \chi \int_0^\infty \io u_k \nabla w_k \cdot\nabla \varphi
	+ \int_0^\infty \io u_k w_k \varphi,
  \ee
  by Definition \ref{defi_weak} known to be valid for all $k\in\N$, we may choose $k=k_l$ and take $l\to\infty$
  in the first, third and fifth summand to obtain that
  \bas
	- \int_0^\infty \io u_{k_l} \varphi_t \to -\int_0^\infty \io \hatu \varphi_t,
	\qquad \mbox{and} \qquad
	- D \int_0^\infty \io \nabla u_{k_l} \cdot \nabla \varphi
	\to - D \int_0^\infty \io \nabla \hatu \cdot \nabla \varphi
  \eas
  as well as
  \bas
	\int_0^\infty \io u_{k_l} w_{k_l} \varphi
	\to 0
  \eas
  as $l\to\infty$.
  In order to show decay also of the cross-diffusive contribution to (\ref{16.66}), relying on our restriction that
  $\frac{\partial\varphi}{\partial\nu}=0$ on $\pO\times (0,\infty)$
  we once more integrate by parts therein to warrant accessibility to (\ref{16.5}) and (\ref{16.6}),
  which namely imply that indeed
  \bas
	\chi \int_0^\infty\io u_{k_l} \nabla w_{k_l} \cdot\nabla \varphi
	= - \chi \int_0^\infty \io w_{k_l} \nabla u_{k_l} \cdot\nabla \varphi
	- \chi \int_0^\infty \io u_{k_l} w_{k_l} \Delta \varphi
	\to 0
	\qquad \mbox{as } l\to\infty.
  \eas
  Therefore, (\ref{16.66}) entails that for any such $\varphi$ we have
  \bas
	- \int_0^\infty \io \hatu \varphi_t 
	-\io u_0 \varphi(\cdot,0)
	= -D\int_0^\infty \io \nabla \hatu \cdot \nabla\varphi,
  \eas
  so that well-known uniqueness arguments for generalized solutions of the heat equation (e.g.~reasonings proceeding
  by duality, as presented in \cite[Proposition 52.13]{quittner_souplet}) become applicable so as to assert that
  actually $\hatu$ must coincide with $U$.
  The claimed approximation property (\ref{16.4}) therefore results from (\ref{16.55}).
\qed
In conjunction with Lemma \ref{lem15}, this shows that as $w_0$ becomes small in the above flavor, 
the quantity $\io \ln u(\cdot ,t) -\io \ln u_0$ indeed remains uniformly positive within a suitable set of times:
\begin{lem}\label{lem17}
  Assume (\ref{iuv}), and let $L>0$ and $t_0>0$ be as in Lemma \ref{lem15}. 
  Then for all $M>0$ there exists $\sigma_0(M)>0$ such that whenever (\ref{M}) and (\ref{d}) hold with
  $\sigma=\sigma_0(M)$, one can find a measurable set $S\subset (t_0,t_0+1)$ such that $|S| \ge \frac{1}{2}$
  and $\ln u(\cdot,t)\in L^1(\Omega)$ for all $t\in S$ as well as
  \be{17.1}
	\io \ln u(\cdot,t) - \io \ln u_0 \ge \frac{L}{2}
	\qquad \mbox{for all } t\in S.
  \ee
\end{lem}
\proof
  In view of (\ref{124.01}), we see that if the claim was false then there would exist a sequence
  $(w_{0k})_{k\in\N} \subset W^{1,\infty}(\Omega)$ of nonnegative functions $w_{0k}$ fulfilling
  $\sqrt{w_{0k}} \in W^{1,2}(\Omega)$ and
  \be{17.2}
	\io \frac{|\nabla w_{0k}|^2}{w_{0k}} \le M
	\qquad \mbox{for all } k\in\N
  \ee
  as well as
  \be{17.3}
	w_{0k} \to 0
	\quad \mbox{in } L^\infty(\Omega)
	\qquad \mbox{as } k\to\infty,
  \ee
  but such that for all $k\in\N$, the measurable set
  \bas
	S_k:= \bigg\{ t\in (t_0,t_0+1) \ \bigg| \
	\ln u(\cdot,t) \in L^1(\Omega) \mbox{ with } \io \ln u_k(\cdot,t) - \io \ln u_0 <\frac{L}{2} \bigg\}
  \eas
  would satisfy
  \be{17.4}
	|S_k| > \frac{1}{2},
  \ee
  where $(u_k,v_k,w_k)$ denotes the associated global weak solution of (\ref{0}) from Theorem \ref{theo222}
  with $w_0:=w_{0k}$.
  Now passing to a subsequence if necessary, from Lemma \ref{lem16} would obtain a null set $N\subset (0,\infty)$
  such that (\ref{16.3}) holds, and such that
  \bas
	f_k(t):=\io \ln u_k(\cdot,t) - \io \ln u_0,
	\qquad t\in (0,\infty)\setminus N, \ k\in\N,
  \eas
  satisfies
  \bas
	f_k(t) \to f(t):=\io \ln u(\cdot,t) - \io \ln u_0
	\quad \mbox{for all } t\in (0,\infty)\setminus N
	\qquad \mbox{as } k\to\infty.
  \eas
  Since thus $f_k\to f$ a.e.~in the bounded interval $(t_0,t_0+1)$ as $k\to\infty$, Egorov's theorem would apply 
  so as to assert that $f_k \to f$ almost uniformly in this interval, in particular meaning that we could find $k_0\in\N$
  and a measurable set $E\subset (t_0,t_0+1) \setminus N$ such that $|E| \ge \frac{3}{4}$ and
  \bas
	\Big| f_{k_0}(t) - f(t) \Big| \le \frac{L}{4}
	\qquad \mbox{for all } t\in E.
  \eas
  As thanks to Lemma \ref{lem15} our choices of $L$ and $t_0$ warrant that
  \bas
	f(t) \ge L
	\qquad \mbox{for all } t\ge t_0,
  \eas
  this entails that
  \bas
	f_{k_0}(t) \ge \frac{3L}{4}
	\qquad \mbox{for all } t\in E
  \eas
  and that hence $E\cap S_{k_0}=\emptyset$.
  This, however, is possible only if $|E| + |S_{k_0}| \le |(t_0,t_0+1)|=1$ and thus requires that
  $|S_{k_0}| \le 1-|E| \le \frac{1}{4}$ which is inconsistent with (\ref{17.4}) and thereby proves the lemma.
\qed
Combining this with Lemma \ref{lem14}, we finally arrive at our main result on migration-driven advantage
of the first population in comparison to the static one:\abs
\proofc of Theorem \ref{theo18}. \quad
  We let $\sigma_0(M)>0$ be as given by Lemma \ref{lem17}, and taking $a>0, b>0, L>0$ and $t_0>0$ from (\ref{ab})
  and Lemma \ref{lem15}, respectively, we fix $\sigma>0$ small enough fulfilling $\sigma\le \sigma_0(M)$ and
  \bas
	a|\Omega|\sigma + \Big( b + \frac{\chi^2}{32 D} \Big) |\Omega| \sigma^2 \le \frac{L}{16}.
  \eas
  Thus, assuming $w_0$ to be compatible with (\ref{iw}), (\ref{M}) and (\ref{d}), from Lemma \ref{lem14} we obtain a 
  null set $N\subset (0,\infty)$ such that the global weak solution of (\ref{0}) from Theorem \ref{theo222} satisfies
  (\ref{14.01}) and
  \bea{18.3}
	\io \ln \frac{v(\cdot,t)}{u(\cdot,t)}
	&\le& \io \ln \frac{v_0}{u_0}
	+ a\io w_0 + \Big(b+\frac{\chi^2}{32D}\Big) \io w_0^2
	- \frac{1}{4} \cdot \bigg\{ \io \ln u(\cdot,t_\star) - \io \ln u_0 \bigg\} \nn\\
	&\le& \io \ln \frac{v_0}{u_0}
	+ \frac{L}{16} 
	- \frac{1}{4} \cdot \bigg\{ \io \ln u(\cdot,t_\star) - \io \ln u_0 \bigg\} \nn\\[2mm]
	& &  \hspace*{30mm}
	\qquad \mbox{for all $t_\star \in (0,\infty)\setminus N$ and any } t\in (t_\star,\infty)\setminus N.
  \eea
  On the other hand, since $\sigma\le \sigma_0(M)$ we may invoke Lemma \ref{lem17} to obtain a measurable set
  $S\subset (t_0,t_0+1)\setminus N$ such that $|S|\ge \frac{1}{2}$ and that (\ref{17.1}) is valid, in particular
  ensuring the existence of $t_\star\in (t_0,t_0+1)\setminus N$ such that
  \bas
	\frac{L}{16}
	- \frac{1}{4} \cdot \bigg\{ \io \ln u(\cdot,t_\star) - \io \ln u_0 \bigg\}
	\le \frac{L}{16} - \frac{L}{8} = - \frac{L}{16}.
  \eas
  Therefore, (\ref{18.3}) implies that (\ref{18.1}) holds if we let $C:=\frac{L}{16}$ and $T:=t_0+1$, whereby
  the proof is completed.
\qed

\mysection{Numerical experiments}

Below we provide results of several numerical simulations. Our aim is to get insight into behaviour of the solution to  system \eqref{0}, discussed in Theorem~\ref{theo18}. In particular, we are interested in the evolution in time of $I(t)$ from (\ref{18.1})
%$$
%I(t) = \int_{\Omega} \log \frac{v(x,t)}{u(x,t)}\, dx
%$$ 
in the case which --- without spatial movement --- would favour species $v$, i.e. for $\alpha > \delta$. We set $D_w = \delta=1$,  $\alpha=2$,  $\beta =\gamma = 200$  %several parameters, see Table~\ref{tab:1}, 
and manipulate only $D_u$, $\chi$, together with the initial shapes of $u_0,v_0,w_0>0$, observing the values of 
$$
M_*=\io \frac{|\nabla w_0|^2}{w_0} \qquad \text{ and } \qquad \sigma_* = \|w_0\|_{L^\infty(\Omega)},
$$
i.e. the smallest possible constants in~\eqref{M} and \eqref{d}~in Theorem~\ref{theo18}. We always choose the initial conditions in such a way that $I(0)=0$. 

For numerical solution,  \eqref{0} is discretized with a second order BDF scheme in time. On each time step we discretize the equations with third degree finite elements in space, using the FEniCS package~\cite{fenics-book}.

\begin{figure}
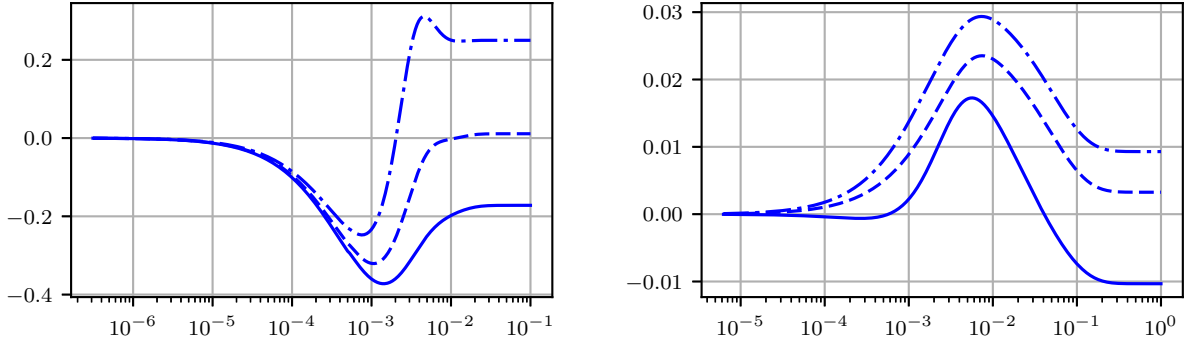
%[ht]
	\centering
	\begin{minipage}{1.2\columnwidth}
		\input{wwk-N-550-S-1600-fig1.pgf}	
		\input{wwk-N-550-S-1600-logpopfig2.pgf}	
	\end{minipage}
	\caption{The evolution of $I(t)$. Left: $w_0(x)=\sigma_*$ with $\sigma_*$ equal to $60$~(solid), $120$~(dashed) or $240$~(dashed--dotted). % After the initial success of $u$ species, either $v$ prevails (dashed--dotted), the initial balance is approximately restored (dashed), or $u$ keeps on the lead (solid). 
		Right:  $w_0(x) = l + (20-l)\exp(-15(x-\frac{1}{2})^2$ with %the same maximum $\sigma=20$ but with varying gradient (and therefore $W_M$). Solid line:~
		$l$ equal to $1.4$~(solid), $14$~(dashed), $20$~(dashed--dotted). See the text for specification of $u_0$ and $v_0$ and other parameters, which are different for both pictures. Note the logarithmic scale of $t$.}
	\label{fig:1}
\end{figure}

\paragraph{Fixed $M_*=0$, varying $\sigma_*$.}

Let us choose $\Omega = [0,1]\subset R^1$ and $\chi = \frac{1}{2}$ and $D_u = 20$. Let us set $u_0 = v_0 = \exp(-15(x-\frac{1}{2})^2)$ as the initial profile\footnote{These initial conditions satisfy the boundary conditions only approximately.}, that is, a Gaussian peak located in the center of~$\Omega$. Assume that $w$ is constant equal to $\sigma_*>0$, so that $M_*=0$. From the left picture in Figure~\ref{fig:1} it turns that depending on the value of $\sigma_*$, after sufficiently long time,  $I(t)$ stabilizes either above or below zero. In particular, when there is not enough food (a situation quite common in nature, when starvation is the typical status of species), the ability to move outside occupied position certainly gives an edge and $I(t)$ remains negative throughout the simulation. On the other hand, when the faster growing species $v$ has enough food, it will finally outgrow $u$ in the sense that $I(t) > 0$ for large~$t$. 

\paragraph{Fixed $\sigma_*$, varying $M_*$.}

Next, for $D_u = 1$ and $\chi = \frac{1}{2}$ we experiment with $u_0$ and $v_0$ with peaks at the opposite ends of $\Omega=[0,1]$,
$
u_0(x) = 1+\exp(-15x^2), \quad v_0(x) = 1+\exp(-15(x-1)^2),
$
and the peak of food distribution $w_0$ located in the middle:
$
w_0(x) = l + (20-l)\exp(-15(x-\frac{1}{2})^2),
$
with $l \in [0,20]$. Note that $\sigma_* = 20$ regardless of $l$, while $l$ obviously affects the value of $M_*$. The right graph in Figure~\ref{fig:1} shows the behaviour of $I(t)$ for $l = 1.4, 14, 20$, with the corresponding values of $M_*$ approximately equal to 170, 10, 0. Note that in this setting, smaller $M_*$ gives bigger advantage to the non-moving species $v$.

Let us remark that for prescribed $M_*$ and $\sigma_*$, various types of behaviour of $I(t)$ are possible, depending on the shape of the initial profiles, so these two parameters alone cannot predict the evolution of the system.

\paragraph{Radially symmetric solution.}

\begin{figure}
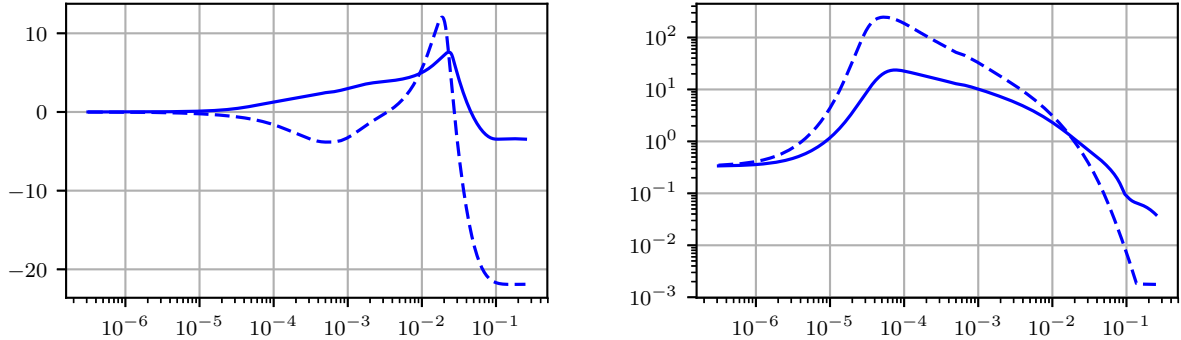
%[ht]
	%	\myplot{0.3}{wwk-test33y--BDF2-p1d-Du-1-Dw-0_05-beta-10-gamma-10-delta-1-alpha-2-chi-1000-T-1--logpop}
	\centering
	%	\myplot{0.3}{wwk-test33y--BDF2-Du-1-Dw-0_05-beta-10-gamma-10-delta-1-alpha-2-chi-1000-T-1--logpop}
	\begin{minipage}{1.2\columnwidth}
		\input{wwk-N-550-S-1600-logpopfig3.pgf}	
		\input{wwk-N-550-S-1600-maxgradufig3.pgf}	
	\end{minipage}
	\caption{Behaviour of the solutions for the same initial profiles in $\Omega$ a unit ball in $R^1$ (solid lines) or in $R^3$ (dashed lines). Left: the graph of $I(t)$ with logarithmic scale of $t$. Right: the log-log graph of $t \mapsto \|\frac{\partial u}{\partial r}(r,t)\|_{L^\infty(\Omega)}$.}
	\label{fig:3}
\end{figure}

Finally, for $D_u = 20$ and $\chi = 10^3$ we solve \eqref{0} in a unit ball $\Omega$ in $R^d$, where $d\in\{1,3\}$, assuming radially symmetric initial data: $u_0(r) = v_0(r) = \frac{1}{10}\exp(-15r^2)$ and $w_0(r) = 2\exp(-15r^2)$, where $r$ denotes the distance from the origin. It follows from Figure~\ref{fig:3} that the solution in $R^3$ admits  larger gradients than the corresponding solution in $R^1$. While in both cases $I(t)$ stabilizes on the negative side, it attains a lower level in $R^3$, leading to an intuition that in ``larger'' space the $u$ species has more space to move to, and to grow.

\mysection{Appendix}
\subsection{Migration-free dynamics: Decisive role of reproduction  rates for ODE asymptotics}
In order to substantiate our discussion around (\ref{0win}), let us briefly address the ODE system related
to (\ref{0A}), as given by
\be{0ODE}
    	\left\{ \begin{array}{rcll}
	u_t &=& \delta uw,
	\qquad & \ t>0, \\[1mm]
	v_t &=& \alpha vw,
	\qquad & \ t>0, \\[1mm]
	w_t &=&  - \beta uw - \gamma vw,
	\qquad & \ t>0, \\[1mm]	
	\end{array} \right.
\ee
with positive parameters $\alpha, \beta, \gamma$ and $\delta$, under
initial conditions $u(0)=u_0\,,v(0)=v_0$ and $w(0)=w_0$.
\begin{prop}
  If $u_0, v_0$ and $w_0$ are positive, then (\ref{0ODE}) admits a unique global positive solution $(u,v,w)$,
  and there exist nonnegative numbers $u_\infty$ and $v_\infty$ such that
  \be{limw} 
	u(t) \to u_\infty,
	\quad
	v(t) \to v_\infty
	\quad \mbox{and} \quad
	w(t) \to 0
	\qquad \mbox{as } t\to\infty.
  \ee
  Moreover, if $v_0=u_0$ then (\ref{0win}) holds, that is, 
  ${\rm{sgn}}(u_\infty-v_\infty) ={\rm{sgn}}(\delta -\alpha)$.
\end{prop}
\proof 
  The statements on existence, uniqueness and positivity are obvious thanks to the Picard-Lindel\"of theorem and a
  comparison argument, because the right-hand sides of (\ref{0ODE}) are locally 
  Lipschitz continuous with respect to $(u,v,w)$, and because, as can easily be seen, the identity
  $ \frac{\beta}{\delta}u +\frac{\gamma}{\alpha} v +w \equiv \frac{\beta}{\delta}u_0 +\frac{\gamma}{\alpha} v_0 +w_0$
%  \be{B1}
%	\frac{\beta}{\delta}u +\frac{\gamma}{\alpha} v +w \equiv 
%	\frac{\beta}{\delta}u_0 +\frac{\gamma}{\alpha} v_0 +w_0
%  \ee
  holds as long as the solution exists.
 Since by positivity it is clear from (\ref{0ODE}) that $u,v$ and $-w$ are nondecreasing we deduce (\ref{limw}).
  Here in the particular case when $u_0=v_0$, on separately integrating the first two equations in (\ref{0ODE})
  we obtain that
  $\ln \frac{v(t)}{u(t)} = {(\alpha- \delta)} \int_0^t w(s)ds \ge (\alpha -\delta) \int_0^1 w(s)ds$ for all $t\ge 1$,
%  \bas
%	 \ln \frac{v(t)}{u(t)} = (\alpha- \delta) \int_0^t w(s)ds
%	\ge (\alpha -\delta) \int_0^1 w(s)ds
%	\qquad \mbox{for all } t\ge 1,
%  \eas
  which directly shows that $u_\infty$ and $v_\infty$ have the claimed ordering property.
%\end{proof}
%
\qed
%
%
%
%
%\begin{rem}
%It follows from the  proposition that in the case of equal  initial conditions for two competing populations  the one with 
%higher  proliferation coefficient   becomes more numerous then the other despite of the values of consumption rates. 
%\end{rem}
%
%
%
%
%
%
%
%
\subsection{A strict form of Jensen's inequality}
Since we could not find an appropriate reference, let us include a brief argument 
for the following essentially well-known result in which we abbreviate $\overline{\varphi}:=\frac{1}{|\Omega|} \io \varphi$
for $\varphi\in L^1(\Omega)$.
\begin{prop}\label{prop999}
  Let $J\subset \R$ be an open interval and $\Psi \in C^2(J)$ be strictly concave. Then
%  for all nonconstant $\varphi\in C^0(\bom)$,
  \be{999.1}
	\Psi (\overline{\varphi}) > \overline{\Psi(\varphi)}
	\qquad \mbox{for all nonconstant } \varphi\in C^0(\bom;J).
  \ee
\end{prop}
\proof
  According to the Taylor theorem, for each $x\in\Omega$ one can find $\xi(x)\in J$ such that
  due to the strict concavity of $\Psi$,
  \bas
	\Psi(\varphi(x))
	- \Psi(\overline{\varphi})
	- \Psi'(\overline{\varphi}) \cdot (\varphi(x) - \overline{\varphi})
	= \frac{1}{2} \Psi''(\xi(x)) \cdot (\varphi(x)-\overline{\varphi})^2 
	\le 
	- c_1 (\varphi(x)-\overline{\varphi})^2
  \eas
  with $c_1:=\frac{1}{2} \min_{\varphi(\bom)} |\Psi''|>0$.
  After integration, this yields
  $ \overline{\Psi(\varphi)} - \Psi(\varphi) - 0 \le c_1 \io (\varphi(x)-\overline{\varphi})^2 dx$
%  \bas
%	\overline{\Psi(\varphi)}
%	- \Psi(\varphi)
%	- 0 
%	\le c_1 \io (\varphi(x)-\overline{\varphi})^2 dx
%  \eas
  and thereby implies (\ref{999.1}), because $\varphi\not\equiv const$.
\qed

\vspace*{5mm}

%
%
%
%
%\noindent
%{\bf Acknowledgement.} \quad
  %The first author acknowledges support of the {\em Deutsche Forschungsgemeinschaft} in the context of the project
  %{\em Analysis of chemotactic cross-diffusion in complex frameworks}.
%
%
%
%

\end{document}